\documentclass[hidelinks,onefignum,onetabnum]{siamart220329}

\ifpdf
\hypersetup{
  pdftitle={Gaussian Variational Schemes on Bounded and Unbounded Domains},
  pdfauthor={J. A. Actor, A. Gruber, E. Cyr, and N. Trask}
}
\fi

\usepackage{setspace}
\usepackage[margin=1in]{geometry}
\usepackage{mathrsfs}
\usepackage[T1]{fontenc}
\usepackage{graphicx}
\usepackage{tikz-cd}
\usepackage{bm}
\usepackage{bbm}
\usepackage{mathtools}
\usepackage[shortlabels]{enumitem}
\usepackage{boldtensors}
\usepackage{amsmath}
\usepackage{amssymb}
\usepackage{subcaption}
\usepackage{algorithm}
\usepackage[noend]{algpseudocode}
\usepackage[utf8]{inputenc}
\usepackage[english]{babel}
\usepackage{graphbox}

\algrenewcommand\algorithmicrequire{\textbf{Input:}}
\algrenewcommand\algorithmicensure{\textbf{Output:}}

\makeatletter
\DeclareFontFamily{OMX}{MnSymbolE}{}
\DeclareSymbolFont{MnLargeSymbols}{OMX}{MnSymbolE}{m}{n}
\SetSymbolFont{MnLargeSymbols}{bold}{OMX}{MnSymbolE}{b}{n}
\DeclareFontShape{OMX}{MnSymbolE}{m}{n}{
    <-6>  MnSymbolE5
   <6-7>  MnSymbolE6
   <7-8>  MnSymbolE7
   <8-9>  MnSymbolE8
   <9-10> MnSymbolE9
  <10-12> MnSymbolE10
  <12->   MnSymbolE12
}{}
\DeclareFontShape{OMX}{MnSymbolE}{b}{n}{
    <-6>  MnSymbolE-Bold5
   <6-7>  MnSymbolE-Bold6
   <7-8>  MnSymbolE-Bold7
   <8-9>  MnSymbolE-Bold8
   <9-10> MnSymbolE-Bold9
  <10-12> MnSymbolE-Bold10
  <12->   MnSymbolE-Bold12
}{}

\let\llangle\@undefined
\let\rrangle\@undefined
\DeclareMathDelimiter{\llangle}{\mathopen}%
                     {MnLargeSymbols}{'164}{MnLargeSymbols}{'164}
\DeclareMathDelimiter{\rrangle}{\mathclose}%
                     {MnLargeSymbols}{'171}{MnLargeSymbols}{'171}
\makeatother

\newcommand{\expec}[1]{\mathbb{E}\left[#1\right]}
\newcommand{\lr}[1]{\left(#1\right)}

\newcommand{\nn}[1]{\left|#1\right|}
\newcommand{\nni}[1]{\left\|#1\right\|}
\newcommand{\intIP}[2]{\left\llangle #1,#2\right\rrangle}
\newcommand{\IP}[2]{\left\langle #1,#2\right\rangle}

\newcommand{\tr}{\mathrm{tr}}
\newcommand{\R}{\mathbb{R}}
\newcommand{\Ints}{\mathbb{N}}

\DeclareMathOperator*{\argmin}{\arg\!\min}

\ifpdf
  \DeclareGraphicsExtensions{.eps,.pdf,.png,.jpg}
\else
  \DeclareGraphicsExtensions{.eps}
\fi

\newsiamremark{remark}{Remark}
\newsiamremark{hypothesis}{Hypothesis}
\crefname{hypothesis}{Hypothesis}{Hypotheses}
\newsiamthm{claim}{Claim}

\headers{Gaussian Variational Schemes}{J. A. Actor, A. Gruber, E. C. Cyr, N. Trask}

\title{Gaussian Variational Schemes on Bounded and Unbounded Domains\thanks{Submitted to the editors \today.
\funding{Funding Note here.}}}

\author{Jonas A. Actor\thanks{Center for Computing Research, Sandia National Laboratories, Albuquerque, NM 
  (\email{jaactor@sandia.gov}).}
\and Anthony Gruber\footnotemark[2]
\and Eric C. Cyr\footnotemark[2]
\and Nathaniel Trask\thanks{School for Engineering and Applied Sciences, University of Pennsylvania, Philadelphia, PA}
}

\usepackage{amsopn}

\newcommand{\norm}[1]{\left \lVert #1 \right \rVert}
\begin{document}

\maketitle

% REQUIRED
\begin{abstract}
A machine-learnable variational scheme using Gaussian radial basis functions (GRBFs) is presented and used to approximate linear problems on bounded and unbounded domains.  In contrast to standard mesh-free methods, which use GRBFs to discretize strong-form differential equations, this work exploits the relationship between integrals of GRBFs, their derivatives, and polynomial moments to produce exact quadrature formulae which enable weak-form expressions.  Combined with trainable GRBF means and covariances, this leads to a flexible, generalized Galerkin variational framework which is applied in the infinite-domain setting where the scheme is conforming, as well as the bounded-domain setting where it is not.  Error rates for the proposed GRBF scheme are derived in each case, and examples are presented demonstrating utility of this approach as a surrogate modeling technique.

% .....\ANT{finish this: mention applications to whitney forms and high-dim problems}

% Compatible meshfree finite elements based on machine-learnable Gaussian radial basis functions (GRBFs) are presented and used to approximate linear problems on bounded and unbounded domains.  Exploiting the relationship between integrals of GRBFs, their derivatives, and polynomial moments leads to exact quadrature formulae which sidestep the curse of dimensionality, enabling the efficient calculation of mass and stiffness matrices necessary for solving the finite element systems.  This yields an extension to the data-driven exterior calculus (DDEC) which is not limited combinatorially by the dimension of the ambient space, and it is shown that the proposed approach effectively learns .... in a way which is convergent, .., and scalable.
% This is an example SIAM \LaTeX\ article. This can be used as a
% template for new articles.  Abstracts must be able to stand alone
% and so cannot contain citations to the paper's references,
% equations, etc.  An abstract must consist of a single paragraph and
% be concise. Because of online formatting, abstracts must appear as
% plain as possible. Any equations should be inline.
\end{abstract}

% REQUIRED

% REQUIRED

\section{Introduction}

% \ANT{Structure of paper:  Intro + High-Level Problem Statement, Method-Theory + Convergence-Poisson, Numerical Scheme for Poisson + Generalization for harder problems/expressions, Whitney Forms theory + method, Examples -- Whitney forms (3D prob with v.f. which is curl + grad) and high-dim Poisson (10D).} \\
% \ANT{Add box around equations to implement.} \\

Data-driven approaches are increasingly being leveraged for building surrogate models for solving physics equations; many of such methods create surrogate models that are both learnable and meshfree, in the sense that their degrees of freedom are not bound to a fixed mesh-based discretizations. Such methods are attractive for surrogate modeling, with an underlying hypothesis that an adaptive or learned representation can encode spatial correlations in the desired solution in a more effective manner than traditional methods. This gives rise to a host of scientific machine learning methods that ultimately build efficient representations of solutions.

In this light, machine-learned surrogate models can be viewed as meshfree data-driven surrogates, where machine learning ultimately uncovers an efficient adaptive basis. However, these bases are rarely controlled in terms of their approximation or stability capabilities.
Some of these meshfree ML methods, such as PINNs \cite{raissi2019physics, zeng2024rbfpinn} and DeepONets \cite{lu2021learning}, fail in practice to consistently build solutions capable of solving equations up to machine precision \cite{cuomo2022scientific}, in both low- and high-dimensional settings \cite{wang2022L2, hu2024tackling}. Others scale poorly with dimension \cite{actor2024data} and thus are ill-suited for practical problems. 
Importantly, most machine learning methods have only vague convergence results (at best), and existing theoretical guarantees do not match observational performance of such methods \cite{shin2020convergence}. 

In comparison, mesh adaptivity applied to classical finite element methods can improve method performance relative to a fixed discretization \cite{devore1989optimal}, and maintain the underlying convergence rates of finite element theory. However, meshing becomes expensive for high-dimensional problems, where surrogate models become increasingly useful. In such settings, practitioners often turn to  approximation methods that scale well in high-dimensions, such as scattered data approximation and radial basis functions. These schemes have been used in collocation schemes for solving partial differential equations, e.g. \cite{schaback1970using, franke1998solving}. However, the analysis of such schemes is difficult and does not provide the full range of stability, consistency, and coercivity results as compared to standard (mesh-based) finite element methods. Moreover, adaptive methods based on these techniques have yet to be fully explored, to the authors' knowledge.

In this work, we pose a novel machine learning surrogate modeling approach well-suited for high-dimensional problems, including problems posed on infinite domains, via the use of machine-learnable Gaussian radial basis functions (GRBFs).  This choice retains adaptivity while effectively sidestepping the issue of poor scalability: by posing a numerical discretization in terms of a GRBF basis, analytic quadrature is retained in a way that scales only weakly with dimension.  This is due to the intimate connection between integration over the infinite support of the GRBF basis and expectation with respect to a Gaussian measure (c.f. Section~\ref{sec:method}).  Note that this choice amounts to a variational crime when the domain in question is compact and bounded, since quadrature over the bounded domain is approximated with an expectation over the support of the basis.  Conversely, the consistency error introduced by this choice is provably bounded (in the case of Poisson's equation, c.f. \ref{thm:errorbound}) and decays exponentially with the variance of the GRBFs.  Since stability is guaranteed by standard estimates, this yields a convergent scheme for Poisson's equation which is empirically demonstrated on a variety of benchmark examples.  Particularly, it is shown that machiine-learnable GRBFs produce an adaptive method well-suited for problems in high dimensions and replete with convergence rates.  Specifically, the contributions of this work are:
\begin{itemize}
    \item A novel variational method for numerically approximating PDE solutions  
    %weak forms of linear equatinos 
    based on GRBFs, complete with convergence proofs for Poisson's equation on bounded and unbounded domains;
    \item An efficient quadrature-free numerical scheme for this method, whose computational cost scales essentially independently of dimension;
    \item A data-driven adaptive surrogate modeling approach that recovers solutions to Poisson's equations in a trainable GRBF basis;
    % \item A data-driven adaptive (in)finite element method based on Gaussian RBFs, which supports \ANT{spatial compatibility?}, analytic quadrature, efficient scaling with dimension, and ...
    \item Applications 
    %to linear problems on finite and infinite domains 
    in low and high ambient dimensions, illustrating the utility of the proposed approach;  
    % \item A demonstration of how our approach can be used to build mimetic discretizations in an adaptive Gaussian basis, using an augmented Whitney form approach.
    % \item \ANT{anything else?}
\end{itemize}

% As this introduces consistency error which may arrest the convergence rate, Nitsche's method \cite{nitsche1971variationsprinzip, benzaken2022constructing} is applied to recover convergence in the case that the boundary is treated by penalty \ANT{double check this claim.  I don't think we can do anything about the consistency error coming from the interior quadrature.}  We show through a collection of linear benchmark examples that the proposed method is convergent and scalable,

Before proceeding with our main results, we offer a sketch of our overarching approach, including the questions and obstacles that we address in the rest of the paper to make such a method pertinent to practical problems of interest.

\section{The method at a glance}
\label{sec:sketch}

Our approach poses a machine learning optimization problem whose solution is a trial space basis of GRBFs to be used in a finite element method.  By minimizing a residual coming from a combination of governing equations and/or observed data, we rely on machine learning to find an optimal basis of GRBFs which reliably approximates the solution space of the weak-form PDE to be solved.
% to uncover a parameterized representation of the solution that minimizes a residual (either from the governing equations or to observed data): we rely on machine learning to uncover an optimal set of Gaussian radial basis functions (GRBFs) to serve as a basis to represent our solution. When discretized in this GRBF basis
To accomplish this, we solve an an equality-constrained nonconvex optimization problem for the means and covariances of the GRBF basis, 
% \ANT{check that notation is consistent with Alg 5.1}
\begin{equation}\label{eq:ML}
\begin{split}
    \theta^*=&\argmin_{\theta}\nn{u_h -u_{\mathrm{data}}}^2, \\
    &\,\,\mathrm{s.t.}\,\, L_{\theta}[u_h] = f,
\end{split}
\end{equation}
where $\theta$ represents the collection of means and covariances of the GRBF basis $\{\phi_i\}$, $u_h=\sum u^i\phi_i$ is the solution representative, and $L_\theta$ is a potentially parameter-dependent linear operator discretizing the underlying weak-form PDE.  The optimal parameters $\theta^*$ then provide a surrogate for the approximation space containing $u_h$, which can be used ``online'' with new initial conditions to generate physics-compatible low-fidelity approximations to PDE solutions in the adapted basis $\{\phi_i\}$.
Importantly, we postulate that a good surrogate model for a large-scale numerical system can be represented with comparatively few GRBFs, particularly when the means and covariances of these basis functions are calibrated to match a given set of high-fidelity data.

The key insight of our approach is in that we obtain a tractable form for the equality constraints -- to enforce that the residual in our discrete system is minimal -- when we use GRBFs as our discretization. A Galerkin variational schemes posed in a GRBF basis lead to exact, quadrature-free assembly on infinite domains, which can be extended to inexact, quadrature-free assembly on bounded domains with appropriate error quantification.  Since GRBFs are easily parameterized in terms of their means and covariances, this implies their use as machine-learnable ``finite elements'' in a way that is compatible with advantages of traditional meshfree methods, such as high-order convergence and weak dependence on ambient dimension.  

To see how this is the case, 
% to mitigate a combinatorial dependence on dimension lies in the use of Gaussian RBFs as basis functions, since RBF approximation scales favorably with dimension.  
recall the standard GRBF density function $\phi_i\sim\mathcal{N}\lr{~m_i,~C_i}$, $1\leq i \leq n$, distributed with mean $~m_i\in\mathbb{R}^d$ and covariance $~C_i\in\mathbb{R}^{d\times d}$,
\[\phi_i(~x) = \lr{2\pi}^{-d/2}\lr{\det~C_i}^{-1/2}\exp\lr{-\frac{1}{2}\lr{~x-~m_i}^\intercal~C^{-1}_i\lr{~x-~m_i}}.\]
Importantly, notice that the gradient of this density,
\begin{equation}\label{eq:gaussgrad}
\nabla \phi_i(~x) = -\phi_i(~x)~C^{-1}_i\lr{~x-~m_i} \coloneqq ~p_i(~x)\phi_i(~x),
\end{equation}
reduces to a linear multiple of itself.  
By the chain rule, it follows that derivatives of Gaussian RBFs are \textit{always} polynomial multiples of themselves, and therefore any strategy for integrating polynomial-times-Gaussian products on the solution domain $\Omega$, which are simply polynomial moments under a Gaussian measure when $\Omega=\R^d$, will be sufficient for assembling the mass and stiffness matrices necessary for solving the Galerkin weak-form system in question.  
% This gives an analytic strategy for the integration in Gaussian variational problems which is scalable to high dimensions and amenable to backpropagation-based optimization.  
% common integral quantities like mass matrices and stiffness matrices in a GRBF basis can be computed in a quadrature-free manner as polynomial moments under a Gaussian measure. 
In the case of the scalar Poisson equation posed on the unbounded domain $\Omega = \mathbb{R}^d$, for example, the usual Galerkin method seeks a solution $u_h \in V_h = \text{span}\{\phi_i\}_{i=1,\dots,n} \subset H^1(\Omega)$ such that, for any $v_h \in V_h$,
\begin{align}\label{eq:wfpoisson}
    B\lr{u_h,v_h} \coloneqq \int_{\Omega} \nabla u_h \cdot{\nabla v_h} = \int_{\Omega}fv_h \eqqcolon F\lr{v_h}.
\end{align}
% It is well known that this method will converge to ... at the rate that the approximation space $V_h$ \cite{something}.
% $$ a(~u,~v) = f(~v),$$
% where
% $$a(~u,~v) = \int_{\R^d} \nabla~u:\nabla~v,\qquad  f(~v) = \int_{\R^d} ~f \cdot ~v .$$
While dealing with integration on infinite domains is less than straightforward for traditional methods, the advantage of assuming a trial space basis of GRBFs lies in assembly: for $\phi_i \phi_j = z_{ij} \phi_{ij}$, with $\phi_{ij}$ the resulting density for the product of Gaussians $\phi_i$ and $\phi_j$ and $z_{ij}$ the corresponding normalization constant, it follows that the entries $\intIP{\phi_i}{\phi_j}=\expec{z_{ij}} = z_{ij}$ and $B(\phi_i,\phi_j)= z_{ij}\expec{ ~p_i\cdot~p_j}$ of the mass and stiffness matrices require only computing the expectations of multivariate polynomials under the Gaussian measure described by $\phi_{ij}$ (c.f. Lemma~\ref{lem:gaussprod}).  Since these have
% described by the product $\phi_i \cdot \phi_j$, 
closed-form, 
% quadrature-free, 
analytic expressions which are explicitly independent of the dimension $d$ (c.f. Theorem~\ref{thm:errorbound}), this suggests that, at least in the case of problems with  polynomial terms and infinite support, GRBF-based Galerkin methods can be successful at mitigating the curse of dimensionality and sidestep mesh-dependent discretizations.
% without committing variational crimes 
%and avoid committing variational crimes for problems posed on all of $\R^d$ 
% \cite{strang1972variational}. 

However, most problems of practical interest are posed on bounded domains, so it is necessary that this technique be equipped to handle this case as well.  Continuing with the example of Poisson's equation, the usual Galerkin method for solving the Dirichlet problem
% \begin{align*}
%     \Delta u &= f \quad \mathrm{in}\,\Omega, \\
%     u &= g \quad \mathrm{on}\,\Gamma,
% \end{align*}
\begin{align*}
    \Delta u = f \quad \mathrm{in}\,\,\Omega, \qquad u = g \quad \mathrm{on}\,\,\Gamma,
\end{align*}
on the bounded domain $\Omega$ with $\partial\Omega=\Gamma$ involves searching for a weak solution $u_h\in V_h$ over some subspace of $H^1_g = \{u\in H^1(\Omega)\,|\,u|_{\Gamma}=g$\}, meaning $H^1$ functions on $\Omega$ with trace $u|_{\Gamma}=g$.  This involves solving the same Equations~\ref{eq:wfpoisson} with test functions $v_h\in V_{h,0}$ which vanish on the boundary $\Gamma$.  Since GRBFs have infinite support and do not approximate $H^1_g$, this prohibits the formulation of a truly Galerkin method with a GRBF discretization.
% Consider as a concrete example the variational form of Poisson's equation, expressed as:
% find $u\in H_g^1(\Omega)$, such that for all $v\in V = H^1_0(\Omega)$,
% \begin{align*}
%     B(u,v) = \int_{\Omega}\nabla u \cdot \nabla v = \int_{\Omega} fv = F(v).
% \end{align*}
Yet, an appropriately generalized Galerkin method is possible, and in view of this, the GRBF scheme proposed here acts fundamentally as a penalty method where the approximation space $V_h\subset C^{\infty}_0(\R^d) \subset H^1_0(\R^d)$ is nonconforming, i.e. $V_h \not\subset H^1_g(\Omega)$.  As a result, the boundary condition $u_h|_{\Gamma}=g$ is not strongly enforced, and the problem to solve is instead: find $u_h\in V_h$, such that for all $v_h\in V_h$,
\begin{align*}
    B_h(u_h,v_h) \coloneqq \int_{\R^d} \nabla u_h\cdot\nabla v_h + \gamma\int_{\Gamma} u_hv_h = \int_{\R^d}\pi_hfv_h + \gamma\int_{\Gamma} gv_h \eqqcolon F_h(v_h),
\end{align*}
where $\pi_hf$ is a projection of the forcing function $f$ onto $V_h$ and $\gamma>0$ is an appropriate penalty parameter penalizing deviation from the desired boundary condition.  As in other cases, the point of the penalty term is to restore coercivity of the bilinear form $B_h$ over the infinite domain $\R^d$, which follows from standard arguments (c.f. \cite{benzaken2022constructing}).  While this penalty formulation is nonstandard in the sense that the solution domain is still the ambient space $\mathbb{R}^d$ and not the bounded domain $\Omega$, it is shown in Section~\ref{sec:fem} that the restriction of the numerical solution $u_h\in V_h$ to $\Omega$ converges to the exact solution $u\in H^1_g(\Omega)$ under appropriate assumptions (c.f. Theorem~\ref{thm:errorbound}), illustrating the extension of the Gaussian variational approach to the bounded case as well.

% Letting $T: V_h \to H^1(\Omega)$ denote the (bounded) truncation operator, the numerical solution $u_h$ resulting from this procedure is simply compared to the true solution $u$ through its truncation $Tu_h$.

% Now that the proposed scheme has been presented in the case of forward simulations, it remains to show that analytic quadrature in combination with the simple parameterization of GRBFs also leads to an effective machine learning scheme for the inverse problem of determining shape parameters which best match a set of observational data.  
% This approach can be extended to problems on compact domainsby borrowing from immersed boundary methods \cite{benzaken2022constructing} to enforce boundary conditions consistently using a penalty method \cite{nitsche1971variationsprinzip, babuska1973finite}. 

More concretely, for the motivating Poisson example, the machine learning problem in \eqref{eq:ML} becomes
\begin{align*}
\theta^* = (~m^*, ~C^*) &= \argmin_{~m, ~C} \sum_{k=1}^{N_\text{data}} \left \lvert \sum_i \widehat{u}_i \phi_i(~x_k) - u_k \right \rvert^2 \\
&\quad\mathrm{s.t.}\,\, L \widehat{u} = F,
\end{align*}
where $L = L(~m,~C)$ and $F = F(~m,~C)$ are given via
% \begin{align*}
% A_{ij}(~m,~C) = A_{ij} &= a(\phi_i,\phi_j) - \gamma \int_\Gamma \nabla_n \phi_i \cdot \phi_j + \int_\Gamma \nabla_n \phi_j \cdot \phi_i + \gamma \int_\Gamma \phi_i \cdot \phi_j \\
% b_i(~m,~C) = b_i &= \int_\Omega f \cdot \phi_i + \gamma \int_{\Gamma} g \cdot \phi_i.
% \end{align*}
\begin{align*}
L_{ij}(~m,~C) &= \int_{\mathbb{R}^d}\nabla\phi_i\cdot\nabla\phi_j + \gamma \int_\Gamma \phi_i \cdot \phi_j \\
F_i(~m,~C) &= \int_{\mathbb{R}^d} \pi_hf\phi_i + \gamma \int_{\Gamma} g\phi_i.
\end{align*}
% In these expressions, $\phi_i,\,A,\,b$ all depend on the shape parameters $~m,\,~C$.
At this point, since this optimization problem is nonconvex, one can turn to their favorite optimizer; as with RBF-Nets, one may choose to use a first-order optimizer favored by machine learning methods, such as Adam, stochastic gradient descent, or limited-memory BFGS.

\begin{remark}
    As in any generalized Galerkin method, the $f$ forcing function built into the right-hand side is projected into the basis of the discrete space $\mathcal{V}_h$ through an appropriate mapping $\pi_h$; alternatively, one can represent $f$ globally via a polynomial or spectral approximation, and then use the exact moment expressions (or Gauss-Hermite quadrature) to evaluate the integral exactly.
\end{remark}

Now that the core ideas have been explained, developing the above procedure into a robust numerical method for problems posed on compact domains 
%(in addition to on all of $\R^d$
% , a number of questions must be satisfactorily answered
requires answering a number of technical questions:
% , which we address in the rest of this paper:
~\\
\begin{enumerate}
\item Under what conditions does a GRBF-based variational formulation converge, and at what rate? 
% At what rate does such a scheme converge?
\item How are boundary conditions precisely treated, and how do they effect the convergence of the underlying scheme?
\item 
% Is there an efficient (i.e., quadrature-free) way to compute the associated integrals? 
If the integration domain is infinite while the problem domain is compact, how does the resulting quadrature error affect the scheme?
% What about other expressions regarding the moments of Gaussians and their gradients?
\item How effectively can we build surrogate models (or otherwise expedite the solution of our problem) by calibrating GRBF shape parameters -- means and covariances -- to match observations of similar solutions?
\end{enumerate}~\\
These questions will be addressed in the remainder of this work, which is organized as follows. Some approximation properties of GRBFs are reviewed in \Cref{sec:approx}, which are used alongside the basic finite element theory reviewed in \Cref{sec:fem} to demonstrate convergence of the proposed method.
% along with a brief description of the necessary finite element theory in Section~\ref{sec:fem}. 
% These notions are used to more precisely formulate the
A more formal treatment of the proposed GRBF variational scheme is given in \Cref{sec:method}, where the necessary quadrature formulae are presented, culminating in Algorithm~\ref{alg:grbf-scheme} which describes the general program.  This program is applied to the scalar Poisson equation in Section~\ref{sec:poissonmethod}, 
%highlighting the benefits and drawbacks of the approach.
% and used to develop a generalized Galerkin method for Poisson's equation.
% discrete weak-form systems.
% Closed-form expressions for the moments of Whitney $k$-forms in an arbitrary Gaussian basis are then derived in \Cref{sec:expressions}. 
% A particular example of how such a numerical method can be used for Whitney form constructions is shown in \Cref{sec:whitney}. 
and experimental results are provided in \Cref{sec:results} to illustrate the advantages and drawbacks of the proposed approach.  Finally, some concluding remarks and directions for future work are presented in \Cref{sec:conclusion}.
%\ANT{May want Whitney to live as a subsection of \Cref{sec:results}.}

\section{Approximation with Radial Basis Functions}
\label{sec:approx}

Here, several notions from meshfree approximation schemes based on GRBFs related to the present approach are recalled.  The expected convergence rates of Gaussian RBF approximation schemes are particularly outlined, since these will be used to compute the error rates seen Section \ref{sec:fem}.  Connections to universal approximation in neural networks are also drawn. For a more complete accounting of RBF approximation, see e.g. \cite{wendland2004scattered}.

Before we proceed further, we make a brief remark on some notation regarding norms. Continuous, integral-valued norms and inner products are denoted with double-bars, i.e. $\norm{\cdot}$ and $\intIP{\cdot}{\cdot}$. When a specific domain is integrated over, it appears as a subscript, leading to the boundary integral norm $\norm{\cdot}_{\partial \Omega}$. In contrast, discrete Euclidean norms and vector inner products are denoted with single-bars, i.e. $\lvert \cdot \rvert$ and $\IP{\cdot}{\cdot}$. Further specialization of both integral and discrete norms, such as RKHS norms or various $p$-norms, are marked accordingly in both cases.

Let $X  = \{~x_1,\dots,~x_n\} \subset \Omega \subseteq \R^d$ be a (possibly infinite) set of scattered data points. Define $h$ to be the \textit{fill distance} between data points, given as $$h = \sup_{~x \in \Omega} \min_{~x_i \in X} \nn{~x - ~x_i}.$$
For clarity, denote the kernel $K(~x,~y) = \phi(~x - ~y)$.
Define the pre-Hilbert space 
$$H_\Omega = \text{span}\{\phi(\cdot - ~x) \,\,:\,\, ~x \in \Omega\} = \text{span}\{K(\cdot, ~x) \,\,:\,\, ~x \in \Omega\}$$
and an inner product
$$\intIP{f}{g}_{\mathcal{K}} = \intIP{\sum_{i=1}^n c_i K(\cdot,~x_i)}{\sum_{j=1}^n d_j K(\cdot,~x_j)}_\mathcal{K} = \sum_{i=1}^{n} \sum_{j=1}^{n} c_i d_j K(~x_i, ~x_j).$$ 
The space $H_\Omega$ can be completed under the norm induced by this inner product to form the \textit{native space}  $\mathcal{N}_{\phi,\Omega}$ 
associated with a radial basis function $\phi$, which is a reproducing kernel Hilbert space with norm $\norm{\,\cdot\,}_\mathcal{K}$ induced by the inner product $\intIP{\cdot}{\cdot}_\mathcal{K}$.

With this, there are the following error estimates for approximation with Gaussian radial basis functions \cite{demarchi2018lectures, fasshauer2007meshfree}.
\begin{theorem}
Let $f \in \mathcal{N}_{\phi,\Omega}$ and let $P_f$ be its GRBF interpolant on $X$. For any $x \in \Omega$, denote $S_r(~x) = \Omega \cap B_{r}(~x)$. Then, for any $k \in \Ints$ there exist constants $C, r$ independent of $X$, $f$, and $\phi$, such that for any $x \in \Omega$, 
\begin{equation}
\lvert f(~x) - P_f(~x) \rvert \le C h^k \sqrt{C_\phi(~x)} \norm{f}_\mathcal{K},
\end{equation}
where $h$ is the fill distance and 
$$C_\phi(~x) = \max_{\lvert \vec{\beta} \rvert + \lvert \vec{\alpha} \rvert = 2k}\,\, \max_{~y,~z \in S_{rh}(~x) } \,\, \left\lvert \, \partial^{\vec{\alpha}}_{~y} \, \partial^{\vec{\beta}}_{~z} \,\, \phi(~y-~z) \, \right\rvert.$$
\end{theorem}
\begin{corollary}
Since the Gaussian density function $\phi \in C^\infty(\R^d)$ and the above estimate holds $\forall k \in \Ints$, taking the limit as $k \rightarrow \infty$ one can derive the well-known ``quadratic exponential'' error bound: \cite{demarchi2018lectures,schaback1997reconstruction} $$\norm{f - Pf}_\infty = \mathcal{O}\left( \exp\left( -\frac{c}{h^2} \right) \right).$$
\end{corollary}
\begin{remark}
Similar bounds exist for relaxed assumptions and some $f \notin \mathcal{N}_{\phi,\Omega}$, but are much more involved; see \cite{fasshauer2007meshfree, schaback1997reconstruction} and references therein. 
Similarly, similar bounds exist in other norms with slightly different assumptions on $\Omega$ and $X$; see \cite{yoon2001spectral, light1998power} and references therein.
When $\Omega \subset \R^d$ is a bounded subset satisfying an interior cone condition, the sharpest error bound at present is 
\begin{equation} \label{eq:GaussInterpLInf}
\norm{f - Pf}_\infty = \mathcal{O}\left( \exp\left( -\frac{c}{h}\right) \right),
\end{equation}
although it is unclear where this result originated;
% if this result has been concretely proven;
see \cite{schaback1997reconstruction, yoon2001spectral} for related discussions pertaining to the provenance of this bound.
\end{remark}

% \subsection{Approximation with RBF-Nets}

Ultimately, the machine-learned, adaptive scheme 
% presented later 
will consider a variant of RBF approximation where the shape parameters (i.e., mean and covariance) both vary anisotropically with each GRBF basis function, contrasting with the standard case of RBF approximation where the covariances are fixed and uniform. Functions of this form are not covered by the classical RBF approximation results mentioned above. Approximation results of shallow MLPs with- such fully trainable RBF shape parameters are discussed in \cite{cybenko1989approximation, stinchcombe1989universal, que2016back}.
% , and a connection to the so-called ``kernel trick'' and Neural Tangent Kernels is mentioned in \cite{zeng2024rbfpinn}.
For our purposes, it is sufficient to note that the error in such approximation schemes is bounded above by the classical case, since one can take the covariances of each trainable Gaussian to be identical to each other.

\section{The GRBF approach as a generalized Galerkin method} \label{sec:fem}

The next goal is to explain the present approach in the context of generalized Galerkin methods, and derive the corresponding fundamental error bound.  To that end, several notions from classical finite elements are now presented which are necessary for accomplishing this task.  More details can be found in, e.g., \cite{brenner2008mathematical, evans2022partial,arnold2010finite}.

% We first in Section \ref{sec:fem} formulate a general framework for Gaussian RBF variational problems, and how such a formulation changes due to boundary conditions in Section \ref{sec:nitsche}. After articulating this theory, we demonstrate how our GRBF scheme is convergent, using the Poisson equation as an example, in Section \ref{sec:poisson}. \ANT{Inconsistent use of subcript $h$ throughout this section.  Also need to make sure function spaces are carefully stated.}

% \subsection{The generalized Galerkin problem}
All problems considered in this work are linear and can be described in the standard weak form, i.e., as seeking a solution $u\in V_g$ in some space of functions $V_g$ satisfying the essential boundary condition $u=g$ on $\Gamma$ which also satisfies the weak-form equations
\[B(u,v) = F(v) \quad \forall v\in V_0,\]
for some bilinear form $B:V\times V\to\mathbb{R}$ and linear functional $F:V\to\mathbb{R}$.  When the bilinear form $B$ is symmetric and positive definite, this can be interpreted as finding the unique Riesz representer for the functional $F$ in the inner product defined by $B$, although any linear problem arising from a variational principle, and many other problems which do not, has a similar weak form.  As mentioned before, choosing a finite dimensional subspace $V_{g,h}\subset V_g$ of the function space $V_g$ and enforcing that both the discretized unknown quantity $u_h\in V_{g,h}$ and the test functions $v_h\in V_{0,h}$ lie there leads directly to the Galerkin problem of finding $u_h\in V_{g,h}$ such that
\[B(u_h,v_h) = F(v_h) \quad \forall v_h\in V_{0,h}. \]
In practice, the space of solutions $V_{h,g}$ is designed to satisfy the essential boundary condition and $V_{0,h}$ is typically comprised of easily parameterizable ``shape functions'', so that the resulting discrete problem can be solved efficiently with a numerical method.  

Conversely, it is not always the case that there are spaces $V_{g,h}\subset V_g$ which are easily parameterizable and/or friendly to numerical computation, particularly when boundaries or boundary conditions $g$ are delicate.  This motivates the use of nonconforming approaches, which rely on finite-dimensional subspaces $V_h\not\subset V_g$.  Since the method presented in Section~\ref{sec:method} uses GRBFs with infinite support as a basis for $V_h$, it will generally be the case that the approximation space $V_h\not\subset V_g$ is not included in the same function space as the exact solution $u$.  This is the essence of a generalized Galerkin scheme, which 
%On the other hand, the method presented in Section~\ref{sec:method} is generalized Galerkin: it 
relies on a sequence of spaces $V_h$ which are not necessarily subspaces of $V_g$, and bilinear resp. linear forms $B_h:V_h\times V_h\to\mathbb{R}$ resp. $F_h:V_h\to\mathbb{R}$ which are not necessarily the restrictions of their continuous counterparts to the discrete space $V_h$.  Assuming continuity of the maps $F_h\mapsto u_h$ (i.e., well-posedness) and defining a formal solution operator $L_h:V_h\to V_h^*$ through $Lu_h(v_h) = B_h(u_h,v_h)$, the generalized Galerkin problem can be formulated as seeking $u_h\in V_h$ satisfying 
\[L_hu_h(v_h) = B_h(u_h,v_h) = F_h(v_h) \quad \forall v_h\in V_h,\]
or, simply as the linear system $L_hu_h = F_h$.  Note that this problem is indeed generalized from the Galerkin problem seen before: choosing a restriction mapping $\pi_h:V\to V_h$, the representative $\pi_hu\in V_h$ of the exact solution $u\in V_g$ may not satisfy the discrete equations, i.e. $L_h\pi_hu\neq F_h$, and the discrete operators $L_h, F_h$ need not be the restrictions of $L,F$.  On the other hand, when the problem is Galerkin and $\pi_h$ is orthogonal projection, the triangle equality yields a straightforward error bound
\begin{align*}
    \nni{u-u_h} &\leq \nni{u-\pi_hu} + \nni{\pi_hu-u_h} \\
    &\leq \nni{u-\pi_hu} + \nni{L_h^{-1}}\nni{L_h\pi_hu-F_h} \\
    &= \inf_{v\in V_h}\nni{u-v} + C_h\nni{L_h\pi_hu-F_h}. 
\end{align*}
Here, $C_h = \nni{L_h^{-1}} = \sup_{F:V\to\mathbb{R}} \nni{u_h}/{\nni{F}}$ is the stability constant of the method while the difference $\nni{L_h\pi_hu-F_h}$ is known as its consistency error and measures how closely the discrete data $L_h,F_h$ approximate their continuous counterparts.  In the Galerkin case, the consistency error can be further estimated,
\[\nni{L_h\pi_hu - u_h} = \sup_{\nni{v_h}=1} \nn{B(\pi_hu,v_h) - F(v_h)} = \sup_{\nni{v_h}=1} \nn{B(\pi_hu-u,v_h)} \leq \nni{B}\inf_{v\in V_h}\nni{u-v},\]
yielding the fundamental error bound for Galerkin methods:
\begin{align*}
    \nni{u-u_h} \leq \lr{1+C_h\nni{B}}\inf_{v\in V_h}\nni{u-v}.
\end{align*}
This shows that the convergence of stable Galerkin methods is entirely controlled by the norm of the bilinear form $B$ and the approximation error of the spaces $V_{g,h}$.  It also applies directly to the proposed GRBF variational method on the unbounded domain $\mathbb{R}^d$, where it is a true Galerkin method.

In the case of generalized Galerkin methods, such as what is needed for a trial space of GRBFs on bounded domains $\Omega$, the situation is more complicated; since $V_h\not\subset V_g$, it is not even obvious what quantity should be compared to the exact solution $u\in V_g$ in order to evaluate the error.  On the other hand, the goal of the method outlined in Section~\ref{sec:sketch} is the approximation on a compact domain $\Omega$ of functions in a Sobolev space $V_g(\Omega)$ of relatively low regularity by functions on an infinite domain $\mathbb{R}^d$ of very high regularity, meaning that $u\in V_g(\Omega) \subset V(\Omega)$ while $V_h\subset C^\infty(\mathbb{R}^d) \subset V(\mathbb{R}^d)$.  Thankfully, there is a bounded extension operator $E: V(\Omega)\to V(\mathbb{R}^d)$ (c.f. Theorem~\ref{thm:extension}) with bounded left inverse $T$ which simply ``ignores'' the function on the complement $\Omega^c$. Therefore, it makes sense to discuss convergence of the quantity
\begin{align*}
    \nni{u-Tu_h}_{V(\Omega)} &\leq \nni{u-T\pi_hu}_{V(\Omega)} + \nni{T\pi_hu-Tu_h}_{V(\Omega)} \\
    &\leq \nni{u-T\pi_hu}_{V(\Omega)} + \nni{\pi_hu-u_h}_{V(\R^d)} \\
    &\leq \nni{u-T\pi_hu}_{V(\Omega)} + \nni{L_h^{-1}}\nni{L_h\pi_hu-F_h}_{V(\R^d)},
\end{align*}
where $\pi_h:V(\Omega)\to V_h$ is, e.g., GRBF interpolation, and the second inequality used that $\nni{T}\leq 1$.  
% Choosing $\pi_h = E$ as the extension operator, this yields the bound
% \begin{align*}
%     \nni{u-Tu_h}_{V(\Omega)} \leq \nni{\pi_hu-u_h}_{V(\mathbb{R}^d)} \leq \nni{L_h^{-1}}\nni{L_hEu - F_h}_{V(\R^d)}.
% \end{align*}
Now, suppose as in Section~\ref{sec:sketch} that $\mathbb{R}^d = \Omega\cup\Omega^c$ and $B_h,F_h$ are chosen in accordance with the penalty method with parameter $\lambda>0$.  It follows that 
\begin{align*}
    B_h\lr{u_h, v_h} &= B\lr{Tu_h,Tv_h} + \lambda \intIP{u_h}{v_h}_{\partial\Omega} + B^c\lr{u_h, v_h} \coloneqq B^\lambda\lr{Tu_h,Tv_h} + B^c\lr{u_h,v_h}, \\
    F_h\lr{v_h} &= F\lr{Tv_h} + \lambda\intIP{g}{v_h}_{\partial\Omega} + F^c\lr{v_h} \coloneqq F^\lambda\lr{Tv_h} + F^c\lr{v_h},
\end{align*}
where $B^\lambda,F^\lambda$ denote the bilinear resp. linear forms on $V(\Omega)$ including the penalty parameter and $B^c,F^c$ denote the remainders of $B_h,F_h$ involving integration over $\Omega^c$.  In particular, notice that $B^\lambda(u,v) = F^\lambda(v)$ when $B(u,v)=F(v)$ and $u\in V_g(\Omega)$. 
% such that $B_h|_{V(\Omega)}=B, F_h|_{V(\Omega)}=F$ \ANT{This isn't quite correct because of the penalty, ugh.  Maybe I can define $B$ to include it.}.  Then, since $\mathbb{R}^d = \Omega\cup\Omega^c$, it follows $\forall v_h\in V_h$ that $B_h(\pi_hu, v_h) = B\lr{T\pi_hu, Tv_h} + B^c\lr{\pi_hu, v_h}$ and  similarly $F_h(v_h) = F(Tv_h) + F^c(v_h)$, where we have defined the (bi)linear functionals $B^c = B_h|_{\Omega^c}, F^c = F_h|_{\Omega^c}$.  
This implies that the consistency error is bounded by
\begin{align*}
    \nni{L_h\pi_hu - F_h}_{V(\R^d)} &= \sup_{v\in V(\mathbb{R}^d)}\frac{\nn{B_h(\pi_hu,v)-F_h(v)}}{\nni{v}} \\
    &\leq \sup_{v\in V(\mathbb{R}^d)}\frac{\nn{B^\lambda\lr{T\pi_hu,Tv}-F^\lambda\lr{Tv}}}{\nni{v}} + \sup_{v\in V(\mathbb{R}^d)}\frac{\nn{B^c\lr{\pi_hu,v}-F^c\lr{v}}}{\nni{v}} \\
    &\leq \nni{B^\lambda}\nni{u-T\pi_hu}_{V(\Omega)} + \varepsilon_{con}\lr{\Omega^c}
\end{align*}
where it was used that $F^\lambda(Tv)=B^\lambda(u,Tv)$ and $\nni{Tv}_{V(\Omega)}\leq \nni{v}_{V(\mathbb{R}^d)}$.  Note that the error
\begin{align*}
    \varepsilon_{con}\lr{\Omega^c} = \sup_{v\in V(\Omega^c)} \frac{\nn{B^c\lr{\pi_hu,v}-F^c\lr{v}}}{\nni{v}},
\end{align*}
measures something akin to consistency error on the complement $\Omega^c$, weighted on this domain by the magnitude of the solution representative $\pi_hu$.  Moreover, the norm of the bilinear form $B^\lambda$ can be bounded in terms of the norm of $B$ with an appropriate trace inequality \cite{brenner2008mathematical}, so that $\nni{B^\lambda} \leq (1+C\lambda)\nni{B}$.  Putting this all together yields and denoting the stability constant as $C_h = \nni{L_h^{-1}}$ as before, the fundamental error bound for our generalized Galerkin method becomes:

\begin{theorem}\label{thm:errorbound}
Suppose $u\in V(\Omega)$ is the unique solution to a well-posed problem of the form: find $u\in V(\Omega)$ such that $B(u,v)=F(v)\,\,\forall v\in V_0(\Omega)$.  Consider the generalized Galerkin problem based on a GRBF discretization: find $u_h\in V_h \subset C_0^\infty(\R^d)$ such that $B_h(u_h,v_h)=F_h(v_h)\,\, \forall v_h\in V_h$, for $B_h,F_h$ defined such that $B_h = B^\lambda + B^c, F_h = F^\lambda + F^c$.  Then, provided this method is stable, the error in the solution on $\Omega$ is bounded as 
\begin{align*}
    \nni{u-Tu_h}_{V(\Omega)} \leq \lr{1+C_h\nni{B^\lambda}}\nni{u-T\pi_hu}_{V(\Omega)} + C_h\varepsilon_{con}\lr{\Omega^c}.
\end{align*}
In particular, the first term is bounded in terms of $\nni{B}$ when there exists a trace inequality $\nni{u}^2_{\partial\Omega}\leq C\,B(u,u).$
\end{theorem}

\begin{lemma}\label{thm:errorbound_rate}
There exist constants $C_1 = C_1(\Omega)>0, c_2>0$ and a restriction mapping $\pi_h$ such that for $u\in N_{\phi,\Omega}$ in the native space associated to $\phi$,
\begin{equation}
\nni{u - T \pi_h u}_{V(\Omega)} \le C_1 \exp\left( - \frac{c_2}{h}\right).
\end{equation}
\end{lemma}
\begin{proof}
With the claim in \eqref{eq:GaussInterpLInf} and identifying the operator $T \pi_h = P$ as the projection operator from Section \ref{sec:approx},
\begin{equation*} \begin{split}
\nni{u - T\pi_h u}_{V(\Omega)} 
&\le \sqrt{\lvert \Omega \rvert} \nni{u - T\pi_h u}_{L^\infty(\Omega)} \le C \lvert \Omega \rvert \exp\left( - \frac{c}{h} \right).
\end{split} \end{equation*}
\end{proof}

The previous results show that the error in the proposed GRBF variational scheme is interpretable and estimable using standard tools from the finite element literature and the approximation theory of GRBFs.  Importantly, note that the global support of the GRBF trial space introduces additional consistency error $\varepsilon_{con}\lr{\Omega^c}$ which does not appear in standard Galerkin methods and must be carefully handled for accurate results.  This will be discussed in detail for the Poisson equation in Section~\ref{sec:poissonmethod}.

\section{Numerical Scheme}\label{sec:method}
It is now possible to discuss the proposed method in more detail. First, note that the use of noncompact GRBFs as elements comes with notable benefits.  Importantly, the choice of GRBFs as basis elements enables exact, quadrature-free integration of polynomial-GRBF products\footnote{We can also exactly integrate trigonometric functions against GRBF densities (see the Supplemental Material), although only the polynomial integration results are needed to build mass and stiffness matrices. } over infinite domains, since \emph{integration} over the infinite domain becomes \emph{expectation} with respect to a product Gaussian measure. To state this more precisely, note the following Lemma proven in Appendix 9, which asserts that the product of a finite number of GRBFs is an (unnormalized) GRBF.

\begin{lemma}\label{lem:gaussprod}
    Suppose $1\leq i \leq n$ and $\phi_i\sim\mathcal{N}\lr{~m_i,~C_i}$,
    \[\phi_i(~x) = \lr{2\pi}^{-d/2}\lr{\det~C_i}^{-1/2}\exp\lr{-\frac{1}{2}\lr{~x-~m_i}^\intercal~C^{-1}_i\lr{~x-~m_i}}.\]
    Then, the product $\prod_{i=1}^{n}\phi_i(~x) = z_{:n}\,\phi_{:n}(~x)\sim z_{:n}\,\mathcal{N}\lr{~m_{:n},~C_{:n}}$ where
    \begin{align*}
        z_{:n} &= \lr{2\pi}^{-\frac{(n-1)d}{2}} \lr{\frac{\det~C_{:n}}{\prod_{i=1}^n \det~C_i}}^{\frac{1}{2}}\exp\lr{-\frac{1}{2}\lr{\sum_{i=1}^n ~m_i^\intercal~C_i^{-1}~m_i - ~m_{:n}^\intercal~C_{:n}^{-1}~m_{:n}}}, \\
        ~C_{:n}^{-1} &= \sum_{i=1}^{n} ~C_i^{-1}, \\
        ~C_{:n}^{-1}~m_{:n} &=\sum_{i=1}^{n}~C_i^{-1}~m_i.
    \end{align*}
    % \begin{align*}
    %     z_{1:n+1} &= \lr{2\pi}^{-d/2} \lr{\det\lr{~C_{1:n} + ~C_{n+1}}}^{-1/2}\exp\lr{-\frac{1}{2}\lr{~m_{1:n}-~m_{n+1}}^\intercal\lr{~C_{1:n}+~C_{n+1}}^{-1}\lr{~m_{1:n}-~m_{n+1}}} \\
    %     ~C_{1:n+1} &= \lr{\sum_{i=1}^{n+1} ~C_i^{-1}}^{-1} \\
    %     ~m_{1:n+1} &= ~C_{1:n+1}\sum_{i=1}^{n+1}~C_i^{-1}~m_i
    % \end{align*}
\end{lemma}

In view of Lemma~\ref{lem:gaussprod}, it follows that the product of GRBF basis functions $\phi_i(~x)\phi_j(~x) = z_{ij}\phi_{ij}(~x)$ becomes a weighted GRBF $z_{ij}\,\mathcal{N}\lr{~m_{ij},~C_{ij}}$, whose weight $z_{ij}$ and parameters $~m_{ij},~C_{ij}$ are fully determined by this result.  As discussed in Section~\ref{sec:sketch}, this means that each entry $0\leq i,j\leq n$ of the scalar-form mass matrix of GRBFs is computable as
\begin{align*}
    \mathcal{M}_{i,j} &= \int_\Omega \phi_i(~x)\phi_j(~x) \approx \int_{\mathbb{R}^d}\phi_i(~x)\phi_j(~x) = z_{ij}\int_{\mathbb{R}^d}\phi_{ij}(~x) = z_{ij}\,\mathbb{E}_{~x\sim\phi_{ij}}[1] = z_{ij},
\end{align*}
which is simply an expectation with respect to the normalized Gaussian measure $\phi_{ij}$ (the case of $\mathcal{M}_{0,0}=\mathrm{Vol}(\Omega)$ is handled exactly).  More generally, it follows from the gradient representation \eqref{eq:gaussgrad} and the chain rule that any product of derivatives (or gradients) of GRBFs is $L^2$ integrable as a polynomial moment in the distribution defined by a product GRBF.  Since it will be shown that these moments are analytically computable, this means that a weak-form discretization based on Gaussian RBFs requires no numerical quadrature, provided some error is accepted in the case of bounded domains $\Omega$ (c.f. Theorem~\ref{thm:errorbound}).

% ; note that this quadrature error decays exponentially with the covariances $~C_i$ and can be made arbitrarily small by adjusting this parameter \ANT{double check this -- erfc function}.  

Generalizing this idea to more Gaussians or higher order polynomial moments required the following definition:
\begin{definition}
    Let $A$ be a tensor with dimensions $[n_1, n_2,...,n_{2p}]$ and let $B$ be a tensor with dimensions $[n_2,n_4,...,n_{2p}]$.  The even-mode tensor contraction $A\times_E B$ is the product with dimensions $[n_1,n_3,...,n_{2p-1}]$ given in components by
    \[\lr{A\times_E B}_{i_1,i_3,...,i_{2p-1}} = \sum_{i_2,i_4,...,i_{2p}}A_{i_1,i_2,...,i_{2p}}B_{i_2,i_4,...,i_{2p}}.\]
\end{definition}

With this in place, it is possible to compute the integral of an arbitrary product of GRBFs and their gradients.  First, let $~{x}^{\otimes d} := ~{x} \otimes ~{x} \otimes \cdots \otimes ~{x}$ denote the outer (tensor) product of $~{x}$ against itself $d$ times. The following proposition shows how to compute the integral of $\alpha$ GRBFs times the outer product of $\beta$ GRBF gradients.

\begin{theorem}\label{thm:quadrature}
Choose (with replacement) $\alpha+\beta$ indices from the set $\{1,\dots,L\}$, and denote them $i_1,\dots,i_\alpha,$ $j_1,\dots,j_\beta$. Let $~X_\ell \sim \mathcal{N}(~{m}_\ell, ~{C}_\ell)$ be Gaussian random variables in $\mathbb{R}^d$.  Let $Z$, $~X$ denote the random variable resp. normalization constant corresponding to the product density $\prod_{a=1}^\alpha\prod_{b=1}^\beta \phi_{i_a}\phi_{j_b}$ (c.f. Lemma~\ref{lem:gaussprod}), and let $~{M}^k = \mathbb{E}(~{X}^{\otimes k})$ denote the $k^{\mathrm{th}}$ tensor moment of $~X$.  Then, it follows that 

\begin{equation}\label{eq:polyquad}
    \begin{split}
        ~{I}^{\alpha,\beta}_{i_1,\dots,i_\alpha;j_1,\dots,j_\beta} &\coloneqq \int_{\mathbb{R}^d} \left( \prod_{a=1}^\alpha \phi_{i_a}(~{x}) \right) \left( \bigotimes_{b=1}^\beta \nabla \phi_{j_b}(~{x}) \right) d~{x} \\
        &= Z\sum_{J \in 2^{\lvert \beta \rvert}} (-1)^{\beta-|J|} \left( \bigotimes_{b=1}^\beta ~{C}_{j_b}^{-1} \right) \times_E \left( \left(\bigotimes_{b=1}^{|J|} ~{m}_{j_b}\right) \otimes ~{M}^{\beta-|J|} \right)_{\sigma_{J}},
        % &= Z \left( \bigotimes_{b=1}^\beta \mathbf{C}_{j_b}^{-1} \right) \times_E \left(\sum_{J \in 2^{\lvert \beta \rvert}} (-1)^{\beta-|J|} \left( \left(\bigotimes_{b=1}^{|J|} \mathbf{m}_{j_b}\right) \otimes \mathbf{M}^{\beta-|J|} \right)_{\sigma_{J}} \right), 
        % &= Z \left( \bigotimes_{b=1}^\beta \mathbf{C}_{j_b}^{-1} \right) \times_E \left(\sum_{J \in 2^{\lvert \beta \rvert}} \sum_{k=|J|}  (-1)^{\beta-k} \left( \bigotimes_{b=1}^k \mathbf{m}_{j_b} \otimes \mathbf{M}^{\beta-k} \right)_{\sigma_{J}} \right), 
        % &= Z \left( \bigotimes_{b=1}^\beta \mathbf{C}_{j_b}^{-1} \right) \times_E \left( \sum_{k=0}^\beta \sum_{\substack{\vec{j} \in 2^{\lvert \beta \rvert} \\ \lvert \vec{j} \rvert = k}} (-1)^{\beta-k} \left( \bigotimes_{b=1}^k \mathbf{m}_{j_b} \otimes \mathbf{M}^{\beta-k} \right)_{\sigma_{\vec{j}}} \right),
    \end{split}
\end{equation}
where $2^{\lvert \beta \rvert}$ is the power set of $\{1,\dots,\beta \}$, $J$ is a multi-index of cardinality $|J|$ representing a subset of $2^{\lvert \beta \rvert}$, and $\sigma_{J}$ is the unique permutation of $\beta$ elements that maps the first $k=\lvert J \rvert$ elements to positions $j_1,\dots,j_k$, and the rest of the elements to positions $j_{k+1},\dots,j_\beta$ \textit{in ascending order}.  Moreover, when $~{y},~{\omega}$ are $N\geq 1$ Gauss-Hermite quadrature nodes and weights, respectively, for polynomials of degree $\beta$, this expression reduces to
\begin{equation}\label{eq:ghquad}
    ~{I}^{\alpha,\beta}_{i_1,\dots,i_\alpha; j_1,\dots,j_\beta} = Z \, \pi^{-\frac{n}{2}}  \sum_{k=1}^N ~{\omega}_k \left( \bigotimes_{b=1}^\beta ~{C}^{-1}_{j_{b}} \left(~{m}_{j_b} - \left( \sqrt{2} ~{C}^{\frac{1}{2}} ~{y}_k + ~{m} \right) \right) \right),
\end{equation}
where $~{C}^{\frac{1}{2}}$ is the Cholesky factor of the covariance matrix $~{C}$ corresponding to $~X$.
\end{theorem}
\begin{proof}
    See Appendix~\ref{app:proofs}.
\end{proof}
\begin{remark}
In the case of using Gauss-Hermite quadrature nodes, only $\beta/2$ quadrature points per dimension are needed for integration over the entire domain; this is in contrast to a pseudospectral or classical FEM quadrature-based approach, which would require this many quadrature points \textit{per element}.
\end{remark}

Theorem~\ref{thm:quadrature} shows how to carry out the necessary integrations forming the core of our GRBF variational scheme.  On the other hand, it was mentioned in Section~\ref{sec:sketch} that boundaries are handled by penalty, meaning that integrals over $\Gamma=\partial\Omega$ are also required.  For this work, these integrals are computed via Monte Carlo sampling; since the dimension of the boundary $\Gamma$ is inherently less than the dimension $d$ of $\mathbb{R}^d$, we are willing to sample sufficiently well along $\Gamma$ to match the error rate of the global method. Moreover, for many high-dimensional applications, including those stemming from the Boltzmann Equation where physical space is bounded but the momentum space is not, the dimension of the boundary may be substantially smaller than that of the full ambient space. Therefore, sampling along the physical 2-D boundary becomes tractable and much more attractive than fully resolving a 6-D solution domain.

Precisely, applying the present GRBF variational scheme relies on the sequence of steps outlined in Algorithm~\ref{alg:grbf-scheme}.

\begin{algorithm}
    \caption{GRBF variational method}
    \label{alg:grbf-scheme}
    \begin{algorithmic}[1]
    \Require A well-posed PDE for the function $u:\Omega\subset\mathbb{R}^d\to\mathbb{R}$ satisfying $u=g$ on $\Gamma=\partial\Omega$.  Data $u_k = u(~x_k)$ ($1\leq k \leq N_{\mathrm{data}}$) for the solution.
    \Ensure An approximation $u_h\in C^\infty\lr{\mathbb{R}^d}$, $u_h = \sum c^i\phi_i$ where $\phi_i = \phi\lr{~m_i,~\Sigma_i}$ is a GRBF with machine-learned parameters $~m_i,~\Sigma_i.$ 
    \State Define a weak-form system  $B(u,v)=F(v)\,\,\forall v\in V_0$ whose solution $u\in V_g(\Omega)$ lies in an appropriate function space.
    \State Define a space of GRBFs $V_h = \mathrm{span}\lr{\phi_i} \subset C^\infty\lr{\mathbb{R}^d}$ and a generalized Galerkin problem on the infinite domain $\mathbb{R}^d$: find $u_h\in V_h$ such that $B_h(u_h, v_h) = F_h(v_h)$ for all $v_h\in V_h$, where the data $B_h = B^\lambda + B^c$ and $F_h = F^\lambda + F^c$ are defined as in Theorem~\ref{thm:errorbound}.
    \State Assemble the discrete problem $L_hu_h = F_h$ according to the quadrature formulae in Theorem \ref{thm:quadrature}.
    \State Letting $~\theta$ denote the collection of means and covariances $\{~m_i,~C_i\}$, calibrate the GRBF basis $\{\phi_i\}$ defining $u_h$ by iteratively solving the problem:
    \[\argmin_{~\theta} \sum_k\nn{u_h(~x_k)-u_k}^2, \quad \mathrm{s.t.}\,\,\,L_{h}u_{h}=F_h,\]
    where $L_h$ and $F_h$ are updated based on the evolving basis $\{\phi_i\}$.

    % \State Solve $L_\theta u_\theta = f_\theta$ and compare the difference $|Tu_{\theta}-u|$.  
\end{algorithmic}
\end{algorithm}
% : \ANT{Want to make this an algorithm or some kind of callout (and clean it up)}
% \begin{enumerate}
%     \item Given a well-posed PDE for the function $u:\Omega\subset\mathbb{R}^d\to\mathbb{R}$ satisfying $u=g$ on $\Gamma=\partial\Omega$, define weak-form equations $B(u,v)=F(v)$ whose solution $u\in V_g(\Omega)$ lies in an appropriate function space.
%     \item Define a generalized Galerkin problem $B_h(u_h, v_h)=F_h(v_h)$ taking place over a space of GRBFs $V_h\subset C^{\infty}(\mathbb{R}^d)$, where the data $B_h = \bar{B}(u_h,v_h)+\lambda\intIP{u_h}{v_h}$ and $F_h = \bar{F}(v_h) + \lambda\intIP{g}{v_h}$ satisfy $\bar{B}|_{\Omega}=B, \bar{F}|_{\Omega}=F$.
%     \item Define a machine-learnable basis $\{\phi_i\}$ spanning $V_h$, and assemble the discrete data $L_\theta, f_\theta$ according to Proposition~\ref{prop:main}.
%     \item Solve $L_\theta u_\theta = f_\theta$ and compare the difference $|Tu_{\theta}-u|$.  
% \end{enumerate}

Once Algorithm~\ref{alg:grbf-scheme} has been carried out for a particular PDE of interest, the resulting basis $\{\phi_i\}$ will be adapted to important characteristics of its solutions.  Therefore, it is reasonable to believe that, after the cost of training is incurred ``offline'', this trained basis can be further useful for predictive simulations ``online''. The remainder of this Section details the application of Algorithm~\ref{alg:grbf-scheme} in the context of the Poisson equation, posed on domains with and without boundary.

% To see how these volumetric quadrature-free integrals work more precisely, we first concretely describe the scheme for Poisson problems in Section \ref{sec:poissonmethod}, and then in Section \ref{sec:poissonML} we articulate the machine learning scheme for building surrogate models using such a formulation. Finally, in Section \ref{sec:arbitrary} we generalize the expressions for integrating products and polynomial moments of Gaussians for arbitrary equations and PDEs discretized in a Gaussian basis. This generalization is necessary for certain applications, such as using the GRBF basis to represent Whitney forms, as done in Section \ref{sec:whitney} later, or for nonlinear problems where the variational problem depends on high-erorder terms.
\section{An Illustration: Poisson's Equation}
\label{sec:poissonmethod}

% \subsection{Example: Convergence in the case of Poisson}
It is now shown that the error bound in Theorem~\ref{thm:errorbound}, along with the approximation results in Section~\ref{sec:approx}, lead to an effective numerical method in the case of Poisson's equation.  First, note the following coercivity property of the relevant bilinear form, which implies stability of the GRBF variational scheme.

\begin{proposition}\label{prop:coercivity}
    The bilinear functional 
    \begin{align*}
        B_h(u,v) = \int_{\mathbb{R}^d} \nabla u\cdot\nabla v + \gamma \int_{\Gamma} uv,
    \end{align*}
    is coercive over $H^1(\Omega)$ and $H^1(\R^d)$.
\end{proposition}
\begin{proof}
    Coercivity on $H^1(\Omega)$ follows from standard arguments.  Using a trace inequality to bound the bilinear form from below (as in \cite{benzaken2022constructing}). Since 
    $$ \norm{ u }_{L^{p}(\partial \Omega)} \le C \norm{ u }_{W^{1,p}(\Omega)} \le C' \norm{ Eu }_{W^{1,p}(\R^d)}, $$
    where $Eu$ is the extension of $u\in W^{1,p}(\Omega)$ (which is guaranteed by Theorem~\ref{thm:extension}),
    it follows that $B_h(u,u) \geq (C')^{-1}\norm{u}^2_{L^2(\Omega)}$ for $u\in H^1(\Omega)$.  Similarly, the case of coercivity over $V_h\subset H^1(\R^d)$ follows immediately from the Gagliardo-Nirenberg-Sobolev
    \footnote{While the statement and proof of the Sobolev inequality in \cite{evans2022partial} is for functions $~u \in C^1_c(\R^n)$, all that is required for the provided proof is that $~u$ be integrable over $\R^n$ along each dimension and that $\lim_{\lvert x_i \rvert \rightarrow \infty} ~u(x) = 0$, which is true of all Gaussian functions.} 
    inequality \cite{evans2022partial}, since
    $$\norm{ u }^2_{L^{2}(\R^d)} \le C \norm{ \nabla u }^2_{L^2(\R^d)} \le B_h(u,u), $$
    implying that $B_h(u,u) \geq C^{-1}\norm{u}^2_{L^2(\R^d)}$ for $u\in H^1(\R^d)$.
\end{proof}

To apply Theorem~\ref{thm:errorbound} in a way which demonstrates convergence, it is necessary to bound the consistency error $\varepsilon_{con}\lr{\Omega^c}$ on the complement of the bounded domain $\Omega$.  A bound on this term will follow from simple Corollaries of Theorem~\ref{thm:quadrature}, which are also proven directly in the Supplemental Material.  
% \ANT{define $\langle \cdot,\cdot \rangle$ meaning here}
\begin{corollary}\label{cor:exxt}
    Let $~x\sim\mathcal{N}\lr{~m, ~C}$, then 
    \[\expec{~x~x^\intercal} = ~C+~m~m^\intercal.\]    
\end{corollary}

\begin{corollary}\label{cor:equadform}
    Let $~x\sim\mathcal{N}\lr{~m, ~C}$, $~a,~b$ be generic vectors, and $~A$ be a generic square matrix.  Then 
    \[\expec{(~x-~b)^\intercal~A(~x-~a)} = \IP{~C}{~A}{} + (~m-~b)^\intercal~A(~m-~a).\]
\end{corollary}

\begin{corollary}
    Let $\phi_i,\,\phi_j$ be GRBF densities related to distributions $\mathbf{X}_i \sim \mathcal{N}(~m_i,~C_i)$ and $\mathbf{X}_j \sim \mathcal{N}(~m_j,~C_j)$, respectively. Then
    \begin{align*}
    \int_{\R^n} \nabla \phi_i \cdot \nabla \phi_j
    &= z_{ij}\,\mathbb{E}\left[ (~x - ~m_i)^\intercal ~C_i^{-1} ~C_j^{-1} (~x - ~m_j) \right] \\
    &= z_{ij}\IP{ ~C_{ij}}{(~C_j ~C_i)^{-1}} + z_{ij}(~m_{ij} - ~m_i)^\intercal ~C_i^{-1} ~C_j^{-1} (~m_{ij} - ~m_j),
    \end{align*}
    where $~C_{ij} = (~C_i^{-1} + ~C_j^{-1})^{-1}$ and $~m_{ij} = ~C_{ij} \left( ~C_i^{-1} ~m_i + ~C_j^{-1} ~m_j \right)$, and the expectation is taken with respect to the normalized density $\phi_{ij} = \phi_i\phi_j/z_{ij}$
\end{corollary} 

With these results in place, the following result bounds the term $\varepsilon_{con}\lr{\Omega^c}$ in Theorem~\ref{thm:errorbound}.

\begin{proposition}\label{prop:consistency}
    Consider GRBFs $\{\phi_i\}_{i=1}^n$ with means $\{~m_i\}$ and covariances $\{~C_i=\sigma^2~I\}$.  Given the approximation space $V_h=\mathrm{span}\{\phi_i\}$ and the operators $B^c: V_h\times V_h\to \mathbb{R}$ and $F^c: V(\Omega)\to\mathbb{R}$, 
    \begin{align*}
        B^c(u,v) = \int_{\Omega^c} \nabla u\cdot\nabla v, \quad F^c(v) = \int_{\Omega^c} Ef,
    \end{align*}
    where $E:V(\Omega)\to V(\mathbb{R}^d)$ is the extension operator from Theorem~\ref{thm:extension}. 
    It follows that the error $\varepsilon^c_{con}(\Omega^c)$ is bounded by
    \begin{align*}
        \varepsilon_{con}^c(\Omega^c) = \sup_{v_h\in V_h}\frac{\nn{B^c(\pi_hu, v_h)-F^c(v_h)}}{\norm{v}} = \mathcal{O}\lr{\sigma^{-(n+4)}\exp\lr{-\frac{c}{2\sigma^2}}}
    \end{align*}
\end{proposition}
\begin{proof}
    We show that the consistency error induced above by integrating over the whole domain decays as the Gaussian variances $\sigma \rightarrow 0$. Let $\bar{f}=Ef$ and note that by the Triangle inequality,
    \begin{align*}
    \left \lvert \int_{\Omega^c}\nabla u_h\cdot\nabla v_h - \int_{\Omega^c} \bar{f} v_h \right \rvert
    &\le \left \lvert \int_{\Omega^c} \nabla u_h\cdot\nabla v_h \right \rvert + \left \lvert \int_{\Omega^c} \bar{f}v_h \right \rvert.
    \end{align*}
    Since Theorem \ref{thm:extension} guarantees control over the smoothness and norm of the extension $\bar{f}$, the second term will decay to zero arbitrarily fast.  In particular, we can choose its decay rate to match that of the first term,
    % The remaining term decays sufficiently quickly as well: for the remaining term,
    \begin{align*}
    \left \lvert \int_{\Omega^c}\nabla u_h
    \cdot\nabla v_h \right \rvert &= \left \lvert \int_{\Omega^c} z_{ij} \left(~C_i^{-1} (~x - ~m_i) \right)^\intercal \left( ~C_j^{-1} (~x - ~m_j) \right) \phi_{ij} \right \rvert.
    \end{align*}
    To determine this rate, assume in the case of uniform refinement that $~C_i = \sigma^2~I$ for all $i=1:n$.
    Define $$~p_{ij}(x) = z_{ij} \left( ~C_i^{-1} (~x - ~m_i) \right)^\intercal \left(~C_j^{-1}(~x - ~m_j) \right),$$
    where $z_{ij},\,m_{ij},\,C_{ij}$ are given by the product-of-Gaussian expressions in Lemma \ref{lem:gaussprod}, and note that 
    \begin{align*}
    ~p_{ij}(~x) &= z_{ij} \sigma^{-4} \left( \nn{~x - ~m_{ij}}^2 - \frac{1}{2}\nn{~m_i - ~m_j}^2 \right).
    \end{align*}
    Thus, the desired integral is bounded by
    \begin{align*}
    \left \lvert \int_{\Omega^c} ~p_{ij}(~x) \phi_{ij}(~x) \right \rvert 
    & \le \int_{\R^d} \left \lvert ~p_{ij}(~x) \right \rvert \phi_{ij}(~x) \\
    &\le z_{ij} \, \sigma^{-4} \int_{\R^d} \left( \nn{~x - ~m_{ij}}^2 + \frac{1}{2}\nn{~m_i - ~m_j}^2 \right) \phi_{ij}(~x) \\
    &= z_{ij} \, \sigma^{-4} \left( \sigma^2 + \frac{1}{2}\nn{~m_i - ~m_j}^2 \right) \\
    &= \mathcal{O}\left( \sigma^{-(n+4)}\exp \left( - \frac{c}{2 \sigma^2} \right) \right),
    \end{align*}
    where we used Lemma \ref{lem:gaussian_moment} to go from the second to the third line.  Note that this expression goes to zero quickly as $\sigma^2 \rightarrow 0$ due to the domination of the exponential term.
    % In particular, we can choose:
    % \begin{align*}
    %     \left \lvert \int_{\Omega^c} \bar{f}v_h \right \rvert
    % \end{align*}
\end{proof}

\begin{remark} We have implicitly assumed here a \textit{unisolvency} requirement that for some $c_1,\,c_2 > 0$, we have related the fill distance $h$ and the GRBF standard deviation $\sigma$ via  $c_1 h \le \sigma \le c_2 h.$  In practice, setting $\sigma = h$ is sufficient for us to see these convergence rates in our computational examples presented later.
\end{remark}

Proposition~\ref{prop:consistency} demonstrates that the additional consistency error contributed by the global support of our GRBFs is small enough not to hinder convergence.  In particular, when $\sigma = c h$ is chosen to be a linear scaling of the fill distance $h$, a convergent error rate is maintained for the overall scheme; as the penalty parameter $\gamma \sim C h \rightarrow \infty$, the rate at which $\gamma$ needs to increase to maintain this rate is determined by traditional trace inequality estimates for this term and match what is seen in practice when the trace coefficient is computed numerically by solving the generalized eigenvalue problem.

Before discussing the applications of this work, it is worth mentioning that the error term $\varepsilon^c_{con}(\Omega^c)$ is complementary to the consistency error introduced by the usual penalty method and cannot be removed with standard Nitsche techniques \cite{evans2022partial,nitsche1971variationsprinzip}.  More precisely, it is easy to see that the standard penalty formulation of the Dirichlet problem on $\Omega$ involving the functionals
\begin{align*}
    \tilde{B}_h(u,v) = \int_\Omega \nabla u\cdot\nabla v + \gamma\int_\Gamma uv, \qquad \tilde{F}_h(v) = \int_\Omega fv + \gamma\int_\Gamma gv,
\end{align*}
is not consistent, since an integration by parts yields
\begin{align*}
    \tilde{B}(u-u_h,v_h) &= \tilde{B}(u,v_h) - \tilde{F}_h(v_h) \\
    &= \int_\Gamma v_h\nabla_{~n}u -\int_\Omega (\Delta u + f)v_h + \gamma\int_\Omega (u-g)v_h \\
    &= \int_\Gamma v_h\nabla_{~n}u \neq 0.
\end{align*}
On the other hand, Nitsche's method provides a technique for modifying the problem data in order to maintain consistency while retaining symmetry in the bilinear form.  In the example above, this looks like:
% This would not be a problem \textit{if we could compute integrals on $\Omega$}, since Nitsche's method defines the new bilinear form and linear functional
\begin{align*}
    B^N(u,v) &= \int_\Omega\nabla u\cdot\nabla v - \int_\Gamma v\nabla_{~n}u + \int_\Gamma u\nabla_{~n}v + \gamma\int_\Gamma uv, \\
    F^N(v) &= \int_\Omega fv + \int_\Gamma g\nabla_{~n}v + \int_\Gamma g v.
\end{align*}
Here, the problem term in the consistency calculation has been added back into the bilinear form to eliminate the consistency penalty, and an additional term has been added to both $B^N,f^N$ in order to symmetrize the bilinear form without changing the weak-form equations.  With these modifications, it is straightforward to check that $B^N(u-u_h,v_h) = 0$, yielding a consistent repair to the penalty method from before.  However, despite the success of Nitsche's method in repairing the problem of finding $u_h$ satisfying $\tilde{B}(u_h,v_h)=\tilde{F}(v_h)$ for all $v_h$, notice that it relies on a variational problem posed on the bounded domain $\Omega$.  This is fundamentally different than the present GRBF approach, which uses the data $B_h,F_h$ containing integrals over $\mathbb{R}^d$.  In this case, there is no additional boundary term that arises in the consistency calculation, since $\Gamma$ has zero measure in the interior of $\mathbb{R}^d$.  However, there is still a contribution to the consistency error as calculated in the discussion above Theorem~\ref{thm:errorbound}, and this is exactly $\varepsilon_{con}\lr{\Omega^c}$.  While this error cannot be eliminated, it is important to note that it will only affect the convergence rate of the scheme, which is of secondary importance for the application of the GRBF variational method since it has been designed as a surrogate modeling technique.  The numerical examples in Section~\ref{sec:results} demonstrate the benefits of this approach that persist despite this peculiarity.

\section{Examples}\label{sec:results}

We now discuss numerical examples which illustrate the proposed GRBF variational approach on bounded and unbounded domains.  To connect to the illustration in Section~\ref{sec:poissonmethod}, the primary exemplar will be Poisson's equation, although it is emphasized that the GRBF variational program in Algorithm~\ref{alg:grbf-scheme} is applicable to any well-posed weak-form PDE of interest (though perhaps with different theoretical guarantees).  Accuracies are reported in units of relative mean squared error (MSE), which is also the loss function used to calibrate the GRBF basis.  Precisely, given space-time data for the true solution $~X$ and an approximate solution $\tilde{~X}$ this metric is:
\[RMSE_{\ell^2} = \frac{\nn{~X-\tilde{~X}}_F^2}{\nn{~X}_F^2},\]
where $\nn{\cdot}_F$ denotes the Frobenius norm.  To facilitate comparison to existing state-of-the-art machine learning methods, the performance of the GRBF variational approach presented here is compared to two different physics-informed neural networks (PINNs) trained to solve the same problems.  The first, denoted PINN (classic) is the basic PINN from \cite{raissi2019physics} which relies on a trial space of fully connected neural networks, while the second, denoted PINN (RBFNet), replaces this trial space with a basis of GRBFs.  Note that the comparisons presented here do not make use of several advanced techniques for training PINNs and other scientific machine learning architectures which are known to improve performance (c.f. \cite{wang2023expert}), 
% such as training with a Neural Tangent Kernel, etc.,
since the goal is to highlight the core approximation capabilities of each network architecture.

The machine-learning implementation is carried out in PyTorch and run on a single Nvidia a100 GPU (41GB).  For each example appearing here, the optimizer used is either Adam or LBFGS (preferred for convergence study), with coefficients in the Gaussian basis constrained as to solve the discretized problem via least squares (as in \cite{cyr2020robust}).  The Adam optimizer runs for 1000 steps with learning rate of $\lambda = 0.01$ and an early stopping criterion that triggers when the loss reaches below $\kappa \cdot 10^{-15}$, where $\kappa$ is the condition number of the discretized linear system.
% We employ a stopping criteria when we reach an error of $\kappa \,\epsilon$, where $\epsilon$ is machine precision and $\kappa$ is the condition number of the linear system. 
The training data in each case is made up of $N=4096$ uniformly distributed points on the domain $\Omega$, and the boundary data used for collocation is also computed by sampling from uniformly-chosen random points on $\Gamma$.  

\subsection{Poisson Problem on Infinite Domain}
Our first test measures convergence, and the benefits of machine learning, on a simple problem: one-dimensional Poisson without boundary. Consider a standard Poisson problem on $\Omega = \R$: find $u: \Omega\to\mathbb{R}$ satisfying  
\begin{align}\label{eq:Poisson}
    -\Delta u &= f, &&\mathrm{in}\,\,\Omega.
\end{align}
    % u &= g &&\mathrm{on}\,\,\Gamma.
Letting $V_h$ denote the GRBF approximation space, the goal is to carry out Algorithm~\ref{alg:grbf-scheme} for $B_h(u_h,v_h) = B(u_h,v_h) = \int_{\mathbb{R}}\nabla u_h\cdot\nabla v_h$ and $F_h(v_h) = F(v_h) = \int_\mathbb{R} fv_h$.  In particular, the GRBF method in this case is a standard Galerkin method.  Choosing $f(x) = \exp \left(-\frac{1}{2} x^2 \right) \lr{x^4 - 2x^3 - 5x^2 + 6x + 2}$, the true solution to this problem is $u(x) = x (2-x) \exp \left(-\frac{1}{2}x^2\right)$. Test data is manufactured by uniformly sampling the interval $\Omega = [-6,6]$, which is taken sufficiently large to capture most of the variance in the solution.

\begin{figure}[htb]
\centering
\includegraphics[width=0.6\textwidth]{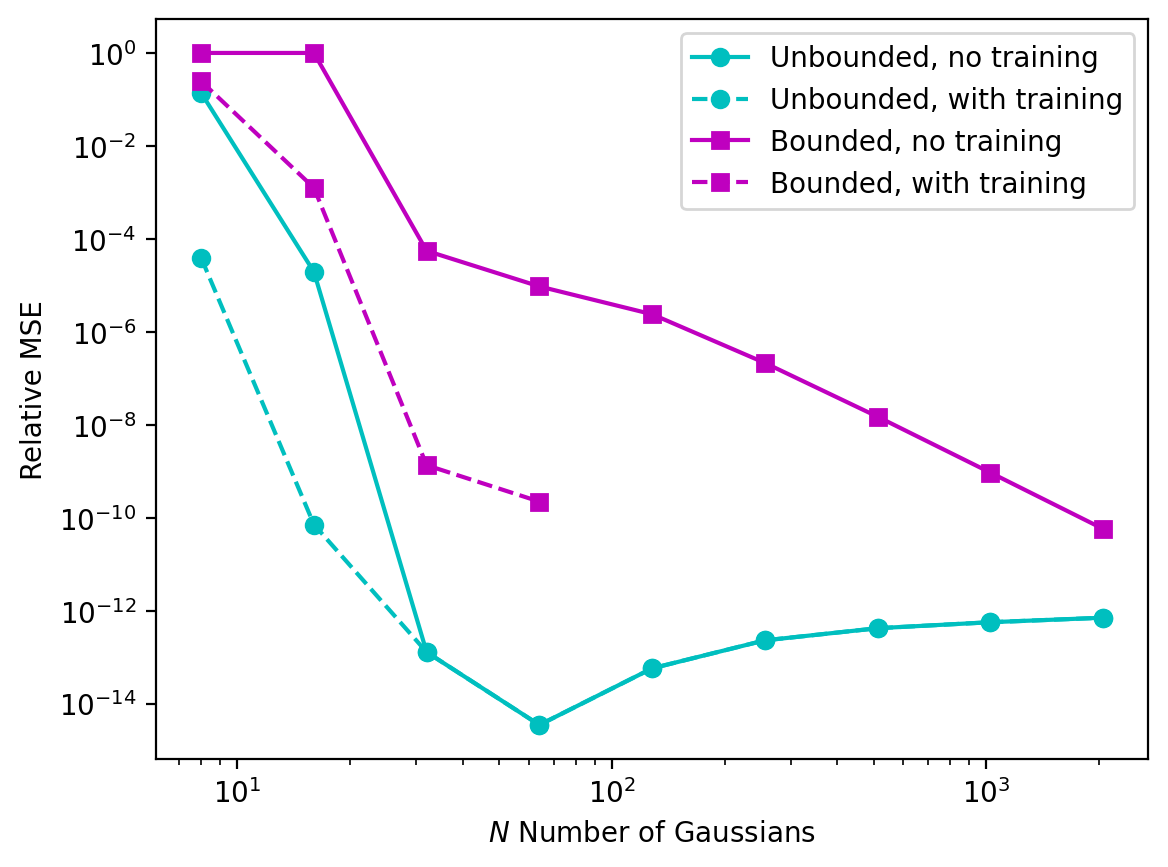}
\caption{Log-log convergence plots for our GRBF method on Problem 1 (unbounded, turquoise) and Problem 2 (bounded, magenta), both with and without training. Note that the error rates of Section~\ref{sec:fem} are observed, although poor conditioning causes a slight climb in the Problem 1 error for large basis sizes. \label{fig:convergence}}
\end{figure}

To test the convergence proven in Section~\ref{sec:fem}, the means of the GRBF basis functions are equispaced uniformly on $\Omega$, with standard deviations scaling as $\sigma= \lvert \Omega \rvert /N$ for $N = 2^k$ with $k=3,\dots,11$. All integrals are computed using the Gauss-Hermite quadrature expressions from Theorem~\ref{thm:quadrature}; all polynomial terms are of polynomial-times-exponential form, and thus require quadrature of degree 3 to integrate exactly. For each $N$, the Gaussian RBF basis is built and directly solved, without training or additional ``$r$-refinement''. Figure~\ref{fig:convergence} shows the relative $\ell_2$ error of the approximation as a function of $N$, with the expected exponential convergence rate observed until machine precision is reached.  Note that the condition number worsens with basis size, causing a slight climb in the error for large numbers of GRBFs.

% Note that this plot is on a log-log scale, demonstrating the error rates proved in Section~\ref{sec:fem}

% \begin{figure}[htb]
% \includegraphics[width=0.45\textwidth]{imgs/convergence_u_loglog}
% \includegraphics[width=0.45\textwidth]{imgs/convergence_F_loglog}
% \caption{Convergence plots for $u$ (left) and $F$ (right), without training to move around Gaussian centers and variances. \label{fig:convergence}}
% \end{figure}

% \subsubsection{Effects of Training}
To measure the benefits of GRBF basis calibration via training, we repeat the problem above, allowing the model to adapt the shape parameters of the GRBFs for 1000 epochs using an LBFGS optimizer and a learning rate of $0.01$. As mentioned at the start of the Section, the accuracy of the GRBF model is compared to a physics-informed neural network (PINN). Note that these approaches are similar in that both build parameterized surrogates of a PDE that scale well to high-dimensional problems. However, there are two inherent differences between the proposed GRBF scheme and PINNs. First, the proposed method solves a variational problem at every step of training and hence guarantees a ``physics residual'' 
% of our discretized system 
of zero by construction.
% (since we directly solve the variational problem at every step of training), 
Conversely, PINNs use the physics residual, i.e. the point a learned operator computed by automatic differentiation on a neural network, as an additional loss term in the optimization problem. Secondly, the GRBF variational method obtains its theorized convergence rates, while PINNs often fail to do so \cite{basir2022critical,rathore2024challenges}. 

To facilitate a fair comparison, two PINN architectures are chosen which are comparably expressive to the GRBF basis, each containing two hidden layers of $N$ neurons. 
% We compare two architectures. 
The first PINN uses a conventional dense MLP with $tanh$ activation function, and has slightly more parameters ($N^2 + (3+d)N$) than the trainable GRBF basis ($dN + \frac{d(d+1)}{2}$). 
The second PINN is a shallow Gaussian RBFNet architecture with $N$ hidden RBF ``neurons'', which possesses a comparable number of parameters to the GRBF model and is closest in structure to the proposed GRBF scheme.
Both PINNs are trained using a combined data fit and physics residual loss, with each term equally weighted.
Otherwise, both PINNs are trained with identical optimizers and hyperparameters as the GRBF scheme above on the same data, using Adam as described at the start of this section.

Table \ref{table:results-prob1} surveys the mean squared error (MSE) from our machine-learnable GRBF method and both of the PINNs, comparing to the static GRBF method from the convergence study.
Since the variational GRBF scheme enforces a physics residual of $0$ by construction, 
% in our discrete space, 
Table \ref{table:results-prob1} reports only the data fit term 
% data-driven component (i.e. not the physics residual term) 
of the MSE from the PINN training.
We see that the MSE quickly converges for the GRBF methods, even without training, reaching machine precision at $N=64$.  Note that training cannot improve upon this accuracy, and there is a slight growth in error for $N>64$ due to the worsening conditioning of the problem.
% for $N \ge 32$, and that training cannot improve upon the machine precision accuracy. 
Conversely, it is seen that the classic PINN does not exhibit any convergence trend as $N$ increases, and that the RBFNet PINN \textit{diverges} as $N$ increases, likely due to increasing numerical instability in the automatic differentiation required to compute the physics residual. Figures showing the resulting models, their accuracy, and the recovered Gaussian bases (for GRBF and the RBFNet PINNs) are provided in the Supplemental Material.
% shown in Section \ref{app:fig} in Figures \ref{fig:prob1-train8-accuracy} and Figure \ref{fig:prob1-train16-accuracy} for $N=8$ and $N=16$, respectively.

\begin{table}[htbp!]
\centering
\begin{tabular}{r|llll}
\hline
& \multicolumn{4}{c}{Model} \\
$N$ & GRBF Solve & GRBF Train & PINN (classic) & PINN (RBFNet) \\
\hline
$8$    & $1.3943 \times 10^{-1}$  & $3.9452 \times 10^{-5}$  & $1.7868 \times 10^{-4}$ & $\mathbf{1.0611 \times 10^{-5}}$ \\
$16$   & $2.0914 \times 10^{-5}$  & $\mathbf{7.2436 \times 10^{-11}}$ & $1.0345 \times 10^{-5}$ & $1.6487 \times 10^{-8}$ \\
$32$   & $\mathbf{1.2705 \times 10^{-13}}$ & $\mathbf{1.2705 \times 10^{-13}}$ & $6.4066 \times 10^{-5}$ & $2.1847 \times 10^{-8}$ \\
$64$   & $\mathbf{3.8740 \times 10^{-15}}$ & $\mathbf{3.8740 \times 10^{-15}}$ & $4.2147 \times 10^{-6}$ & $1.1917 \times 10^{-6}$ \\
$128$  & $\mathbf{5.9538 \times 10^{-14}}$ & $\mathbf{5.9538 \times 10^{-14}}$ & $4.5237 \times 10^{-6}$ & $3.1847 \times 10^{-4}$ \\
$256$  & $\mathbf{2.2211 \times 10^{-13}}$ & $\mathbf{2.2211 \times 10^{-13}}$ & $1.3898 \times 10^{-4}$ & $6.3559 \times 10^{-1}$ \\
$512$  & $\mathbf{4.0239 \times 10^{-13}}$ & $\mathbf{4.0239 \times 10^{-13}}$ & $1.1688 \times 10^{-9}$ & $9.9858 \times 10^{-1}$ \\
$1024$ & $\mathbf{5.5459 \times 10^{-13}}$ & $\mathbf{5.5459 \times 10^{-13}}$ & $1.6252 \times 10^{-5}$ & $1.0021 \times 10^{0}$  \\
$2048$ & $\mathbf{6.2722 \times 10^{-13}}$ & $\mathbf{6.2722 \times 10^{-13}}$ & $3.9468 \times 10^{-1}$ & $4.7978 \times 10^{0}$  \\
\hline
\end{tabular}
\caption{Relative data fit MSE for Problem 1, for an increasing number of Gaussians RBFs / network width, when trained with 1000 steps of Adam. Best recovered MSE for each $N$ is in \textbf{bold}. \label{table:results-prob1}}
\end{table}

\subsection{Poisson Problem on Bounded Domain}

To show that convergence is maintained even in the presence of a boundary penalty 
% the Nitsche penalization of the boundary term, 
we repeat the 1D test case with a manufactured problem posed on the bounded domain $\Omega = [-1,1]$ with $\Gamma = \partial \Omega = \{-1,1\}$.  Choosing a true solution of $u(x) = \sin (3\pi x)$, the right-hand side forcing term $f$ involves trigonometric integrals which admit closed-form expressions similar to the polynomials above; see the Supplemental Material for their derivation. The Gaussian means are initialized with equal spacing on the interval $[-2,2]$, and training data is chosen on the same interval in order to guarantee that $\Gamma$ is well resolved. The same uniform scaling of the Gaussian variances as before is used along with a penalty parameter $\gamma = 16 N$.  In particular, this penalty is chosen which according to the rate at which the largest eigenvalue scales when solving the generalized eigenvalue problem to bound the trace inequality:
\begin{equation}
\int_{\R^n} \nabla \psi_i \cdot \nabla \psi_j = \lambda \int_{\Gamma} \psi_i \, \psi_j,
\end{equation}
with the right-hand side estimated using Monte Carlo sampling. Otherwise, the optimizers and hyperparameters are the same as in the unbounded problem.

Figure \ref{fig:convergence} shows the relative MSE for solving in a Gaussian basis and also the results of training, and Table \ref{table:results-prob2} makes the same comparisons as before between the GRBF methods and the PINNs.  The effects of including boundary terms on the convergence rate are visible in Figure \ref{fig:convergence}, showing the achievement of second-order convergence on this problem. For this set of equations, the GRBF method is difficult to train for large $N$ due to increasingly poor conditioning in the linear solve, while both the standard and RBFNet PINNs fail to train well regardless of $N$. We conjecture that the PINNs' poor performance is due to frequency/spectral bias towards low-frequency content \cite{wang2022NTK}, although additional work would be needed to validate this statement.

% \begin{figure}[htb]
% \centering
% \includegraphics[width=0.6\textwidth]{imgs/sine/sin_convergence_MSE}
% \caption{Convergence plots for Problem 1, both with and without training, for our GRBF method. Note that this plot is on a log-log scale, demonstrating our theoretic error rate. \label{fig:prob2-convergence}}
% \end{figure}

\begin{table}[htbp!]
\centering
\begin{tabular}{r|llll}
\hline
& \multicolumn{4}{c}{Model} \\
$N$ & GRBF Solve & GRBF Train & PINN (classic) & PINN (RBFNet) \\
\hline
$8$    & $1.0000 \times 10^{0}$   & $\mathbf{2.4595 \times 10^{-1}}$  & $9.7701 \times 10^{-1}$ & $9.6239 \times 10^{-1}$ \\
$16$   & $9.9957 \times 10^{-1}$  & $\mathbf{1.2571 \times 10^{-3}}$  & $5.8670 \times 10^{-1}$ & $1.6456 \times 10^{-1}$ \\
$32$   & $5.5429 \times 10^{-5}$  & $\mathbf{1.3605 \times 10^{-9}}$  & $5.1496 \times 10^{-1}$ & $4.9577 \times 10^{-1}$ \\
$64$   & $9.5132 \times 10^{-6}$  & $\mathbf{2.2412 \times 10^{-10}}$ & $3.3363 \times 10^{-1}$ & $1.8007 \times 10^{-1}$ \\
$128$  & $\mathbf{2.3975 \times 10^{-6}}$  & $4.0523 \times 10^{-6} \,\,\,\, \dagger$ & $2.7820 \times 10^{-1}$ & $2.7212 \times 10^{-1}$ \\
$256$  & $\mathbf{2.1560 \times 10^{-7}}$  & $8.4620 \times 10^{-3} \,\,\,\, \dagger$ & $1.2183 \times 10^{-1}$ & $1.1004 \times 10^{-1}$ \\
$512$  & $\mathbf{1.5020 \times 10^{-8}}$  & $7.0509 \times 10^{-1} \,\,\,\, \dagger$ & $1.4696 \times 10^{-1}$ & $7.3170 \times 10^{-1}$ \\
$1024$ & $\mathbf{9.4424 \times 10^{-10}}$ & $3.7859 \times 10^{-1} \,\,\,\, \dagger$ & $6.4065 \times 10^{-1}$ & $1.0066 \times 10^{0}$  \\
$2048$ & $\mathbf{5.8227 \times 10^{-11}}$ & $5.8227 \times 10^{-11}    \,\, \dagger$ & $9.1331 \times 10^{-1}$ & $1.0071 \times 10^{0}$  \\
\hline
\end{tabular}
\caption{Relative data fit MSE for Problem 2, for an increasing number of Gaussians RBFs / network width, when trained with 1000 steps of Adam. Best recovered MSE for each $N$ is in \textbf{bold}. Training runs with GRBFs that reached a stopping criterion due to condition number growth are marked by a dagger $\dagger$. \label{table:results-prob2}}
\end{table}

\subsection{Trained Solutions to High-Dimensional PDEs}

To demonstrate the ability of Algorithm~\ref{alg:grbf-scheme} to handle high-dimensional problems, the next experiment solves a similar
% setup to Problem 1 on a 
Poisson problem posed on an unbounded domain in $\R^d$ for $d=8$. Consider a forcing function of 
$$f(~x) = \exp\left(-\frac{1}{2}\nn{~x}^2\right) \left( 2d + 4\sum_{i=1}^d x_i - 4 \nn{~x}^2 - \left( 2 \sum_{i=1}^d x_i - \nn{~x}^2 \right) \left( \nn{~x}^2 - d \right) \right),$$
corresponding to the manufactured solution 
$$u(~x) = \exp\left(-\frac{1}{2}\nn{~x}^2 \right) \left( \sum_{i=1}^d x_i (2-x_i) \right).$$
In contrast to Problems 1-2, here $65536 = 2^{16}$ points of manufactured data are used (due to the higher latent dimension) along with a learning rate of $\lambda = 0.05$, 10K steps of Adam, and the same ill-conditioning stopping criteria described at the beginning of the Section. As before, Table \ref{table:results-highdim} highlights improvement in the performance of our trained model, from $N=8$ to $N=16$ and then again to $N=32$.  In particular, the ``$r$-refined'' model is capable of relatively low errors with a very small basis size, which is not possible for the untrained model.

\begin{table}[htbp!]
\centering
\begin{tabular}{r|ll}
\hline
& \multicolumn{2}{c}{Model} \\
$N$ & GRBF Solve & GRBF Train  \\
\hline
8 & $9.9406 \times 10^{-1}$ & $1.5789 \times 10^{-5}$ \\
16 & $9.8994 \times 10^{-1}$ & $1.5078 \times 10^{-5}$  \\
32 & $1.0474 \times 10^{0}$ & $7.1046 \times 10^{-6}$  \\
\hline
\end{tabular}
\caption{Relative data fit MSE for Problem 3, for an increasing number of Gaussian RBFs. \label{table:results-highdim}}
\end{table}

\subsection{Whitney Forms}

The final example presented here involves the use of Whitney $k$-forms, which are finite-dimensional differential $k$-forms on a simplicial complex \cite{lohi2021whitney} useful in building compatible discretizations and faithful representations of the de Rham complex corresponding to a PDE system \cite{trask2020conservative, bochev2007principles}. Building surrogate models using Whitney forms is of particular interest due to difficulties with the maintenance of underlying structure, such as geometric invariants, in traditional surrogate modeling techniques \cite{berry2020spectral, gruber2023energetically, actor2024data}.  The GRBF methods developed here are particularly attractive for this task, since other machine-learning approaches for learning Whitney forms, e.g., \cite{actor2024data}, scale poorly with dimension, making them impractical for many real-world applications.

% % traditional surrogate modeling techniques to maintain the geometric invariants and underlying structure of the underlying PDE system \cite{berry2020spectral, gruber2023energetically, actor2024data}. 
% Other machine-learning approaches for learning Whitney forms \cite{actor2024data} scale poorly with dimension, making them impractical for many real-world applications.

% We briefly review key details regarding exterior calculus, Whitney forms, and their role in data-driven surrogates models, before describing how to build Whitney forms in our GRBF scheme.

% \subsection{Exterior Calculus and Whitney Forms}

% \ANT{Double check cochain vs k-form typing}
To describe how Whitney's construction works, recall the notion of a chain complex on the graded vector space $V=\bigoplus_{-\infty}^{\infty}V_k$, which is a sequence of vector subspaces $V_k$ along with linear maps $\partial_k:V_k\to V_{k-1}$ satisfying $\partial_{k+1}\circ\partial_k = 0$, i.e., $\mathrm{im}\,\partial_k \subseteq \mathrm{ker}\,\partial_{k+1}$.  These $V_k$ spaces are commonly chosen as the spaces of $k$-simplices that discretize a domain $\Omega$, and $\partial_k$ is taken to be the boundary operator taking a simplicial chain $v_k\in V_k$ to its boundary chain $\partial v_k\in V_{k-1}$ of one degree lower. These chains can be integrated against continuous differential $k$-forms in such a way that Stokes' theorem holds: for any simplex $v\in V_{k+1}$ and differential $k$-form  $\omega\in \Lambda^k$,
\[\int_{\partial v}\omega = \int_v d_k\omega.\]
Importantly, this defines exterior derivative operators $d_k: \Lambda^k\to\Lambda^{k+1}$ satisfying $d_k\circ d_{k-1}=0$, a property which must be respected by any discretization of differential $k$-forms $\Lambda^k$ in order to maintain discrete analogues of continuous geometric relationships.  Note that $d_0,d_1,d_2$ are the abstract analogues of the usual vector calculus operators $\nabla, \nabla\times, \nabla\cdot$ defined so that $\nabla\times\nabla = \nabla\cdot\nabla\times = 0$.

% The topological duals of the $\partial_k$ under the standard integration pairing, denoted by $\delta_k:V_k^*\to V_{k+1}^*$, are the (discrete) exterior derivatives, which are defined so that Stokes' theorem holds

% Stokes' theorem, expressed 
% for a simplex $v\in V^k$ and a differential $k$-form  $\omega\in\Lambda^k$ as
% \[\int_{\partial v}\omega = \int_v d\omega,\]
% then induces a dual cochain complex in terms of
% the coboundary operators $d_k:V^k\to V_{k+1}$, which are restrictions of the exterior derivative to $k$-forms on $V_k$:

% \ANT{TikZ picture here}

% The spaces $\lambda^k$ can be consistently discretized with Whitney $k$-forms,
In this language, Whitney $k$-forms are a compatible discretization of the continuous spaces $\Lambda^k$ and are defined as follows: for any simplex $v = [p_0,...,p_k] \in V_k$, 
\begin{align*}
    W_v = k! \sum_{i=0}^k (-1)^i \lambda_i \, d\lambda_0\wedge...\wedge\widehat{d\lambda_i}\wedge...\wedge d\lambda_k,
\end{align*}
where $\lambda_i:V_k\to\mathbb{R}$ are the barycentric coordinates and $\widehat{d\lambda_i}$ indicates a term which is left out.  It follows that $d_{k+1}\circ d_k = 0$ for all Whitney forms $\mathcal{W}_v$, as required for de Rham compatibility.  Moreover, the Whitney $k$-forms contain all constant $k$-forms since they form a partition of unity, i.e. $\sum\lambda_i(~x)=1$ for all $~x$.  In the context of machine-learned surrogate models, such barycentric coordinates or partitions of unity can come from a coarsened graph representation of the domain \cite{trask2022enforcing} or otherwise calibrated based on data \cite{actor2024data}, providing a basis which is tailored to the behavior of a particular high-fidelity model. 

% Importantly, the use of Whitney forms ensures that consequences of the cochain complex are exactly preserved, enabling reconstructions which preserve involution constraints and are free of spurious computational modes. In the context of this work, we use GRBFs to define the spaces $V^k$, taking care to preserve the properties of the underlying chain complex.

% There are two essential ingredients: first, there is the use of Gaussian RBFs as Whitney $k$-forms, which leads to the spatial compatibility necessary for an effective mixed method.  Second, there is the treatment of boundaries, which nonobvious due to the use of noncompactly supported RBFs.  

% \subsection{Using Gaussian RBFs as Whitney \texorpdfstring{$k$}-Forms}
% \label{sec:whitney-ex}

% \begin{itemize}
% \item discussion of compatibility of using Gaussians, since no POU property
% \item add $\mathbf{1}$ to basis
% \end{itemize}

While the proposed GRBF basis functions do not form a partition of unity, they can still be used to build Whitney forms.  Defining the characteristic function $\psi_0(~x)=1_{\Omega}$, the set $\{\psi_0,\dots,\psi_n\}$ forms a Whitney 0-form basis for functions on $\mathbb{R}^d$, from which a vector proxy for the Whitney 1-forms can be constructed via
\begin{align*}
    \psi_{ij}(~x) &= \psi_i(~x)\nabla\psi_j(~x)-\psi_j(~x)\nabla\psi_i(~x) = \psi_i(~x)\psi_j(~x)\lr{~p_j(~x)-~p_i(~x)},
    % \psi_{ij}(~x) &= \psi_i(~x)\nabla\psi_j(~x)-\psi_j(~x)\nabla\psi_i(~x) = \psi_i(~x)\psi_j(~x)\lr{~C^{-1}_i(~x-~m_i)-~C^{-1}_j(~x-~m_j)},
\end{align*}
where $0\leq i<j\leq n$ and $~m_0, ~C_0^{-1}$ are treated as identically zero.  Clearly, the anticommutative relationship $\psi_{ji}=-\psi_{ij}$ holds for all pairs of indices $i,j$, as required for a 1-form basis.  Moreover, it follows that the nontrivial gradient of each 0-form is expressible as
\[\nabla\psi_i(~x) = \psi_0(~x)\nabla\psi_1(~x)-\psi_1(~x)\nabla\psi_0(~x) = \psi_{01},\]
so that the gradient of each 0-form is contained in the span of the 1-form basis.  It is straightforward to check that these relationships persists for higher-order $k$-forms, so that GRBF Whitney forms are compatible in the sense of $d_k \circ d_{k-1} = 0$.  Importantly, considering a simplex $v = [p_{i_0},...,p_{i_k}]$ and using the fact that $d\psi_i(~x) = \psi_i(~x)~p_i(~x)\cdot d~x$, it follows that
\begin{equation*}
    W_v = k!\prod_{j=0}^k \psi_{i_j} \sum_{j=0}^k \lr{-1}^j \lr{~p_{i_0}\cdot d~x}\wedge ...\wedge \widehat{~p_{i_j}\cdot d~x}\wedge ... \wedge \lr{~p_{i_k}\cdot d~x},
\end{equation*}
so that the same product of Gaussian PDFs appears in every term of the Whitney $k$-form.  This provides a recipe for constructing GRBF Whitney $k$-forms which will be used to solve the following mixed-form Poisson problem.

% \ANT{probably should state general formula}

\begin{remark}
    The nonstandard inclusion of the constant function $\psi_0=1_{\Omega}$ in the 0-form basis is present because Gaussian RBFs do not form a partition of unity, and therefore cannot represent constants exactly.  Since nonzero constants are not integrable on infinite domains, this means that the present approach to GRBF Whitney forms is defined only for bounded domains $\Omega$.  It is likely that some conceptual modification can be introduced in the future to avoid this restriction.
\end{remark}

Assembly for mass and stiffness matrices 
% for the discrete spaces $~V^k$ using 
corresponding to weak-form expressions in our machine-learned basis GRBF can be done efficiently; we present a few corollaries which illustrate how Theorem~\ref{thm:quadrature} appears in important special cases.  For comparison, direct calculations of these expressions by a different method are given in the Supplemental Material. In the results that follow, we define the $k$-form mass matrix $~M^k$ and stiffness matrix $~S^k$ as
\begin{equation*} \begin{split}
{M}^k_{i_0,\dots,i_k;\,j_0,\dots,j_k} &= \intIP{ \psi_{i_0,\dots,i_k} }{ \psi_{j_0,\dots,j_k}}, \\
{D}^k_{i_0,\dots,i_k;\,j_0,\dots,j_{k+1}} &= \intIP{ d_k \psi_{i_0,\dots,i_k} }{\psi_{j_0,\dots,j_{k+1}}},\\
{S}^k_{i_0,\dots,i_k;\,j_0,\dots,j_{k}} &= \intIP{ d_k \psi_{i_0,\dots,i_k} }{d_k\psi_{j_0,\dots,j_{k+1}}}.
\end{split} \end{equation*}
% For such matrices, we have the following corollaries:
\begin{corollary}
The mass matrix $~M^0$ of 0-forms and stiffness matrix $~S^0$ of 0-form gradients are given component-wise as 
\begin{equation*}
    M^0_{ij} = \intIP{\psi_i}{\psi_j} = I^{2,0}_{ij}, \qquad S^0_{ij} = \intIP{\nabla\psi_i}{\nabla\psi_j} = \tr \left( ~I^{0,2}_{ij} \right).
\end{equation*}
Additionally, the mixed matrix $~D^0$ between 0-form gradients and 1-forms is given component-wise as
\begin{equation*}
    D^0_{i,ab} = \intIP{\nabla\psi_i}{\psi_{ab}} = \tr\lr{~I^{1,2}_{a,bi}-~I^{1,2}_{b,ai}}.
\end{equation*}
\end{corollary}

\begin{corollary}\label{cor:1formM}
The mass matrices $~M^1$ of 1-forms and $~M^2$ of 2-forms are given component-wise as 
\begin{align*}
    M^1_{ij,ab} &= \intIP{\psi_{ij}}{\psi_{ab}} = \tr \left( ~I^{2,2}_{ia,jb} - ~I^{2,2}_{ib,ja} - ~I^{2,2}_{ja,ib} + ~I^{2,2}_{jb,ia} \right), \\
    M^2_{ijk,abc} &= \intIP{\psi_{ijk}}{\psi_{abc}} = \sum_{(ijk)} \sum_{(abc)} \tr_{12} \tr_{34} \left( ~{I}^{2,4}_{ia,jbkc} - ~{I}^{2,4}_{ia,jckb} \right),
\end{align*}
where $\tr_{ij}$ denotes a trace over indices $i,j$ in its tensor argument.  Moreover, the mixed matrix between 1-form curls and 2-forms is given component-wise as
\begin{equation*}
    D^1_{ij,abc} = \intIP{\nabla \times \psi_{ij}}{\psi_{abc}} = 2 \sum_{(abc)} \tr_{12}\tr_{34}\left( ~I^{1,4}_{a,ibjc} - ~I^{1,4}_{a,icjb} \right).
\end{equation*}
\end{corollary}

% \begin{corollary}
% The mixed matrix between 0-forms and 1-forms $\mathbf{V}^0$ is given as $$\mathbf{V}^0_{i,ab} = \langle \nabla \psi_i, \psi_{ab} \rangle = \tr \left( \mathbf{I}^{1,2}_{a,bi}  - \mathbf{I}^{1,2}_{b,ai}\right).$$
% \end{corollary}
% % \color{red}{
% \begin{corollary}
% The mixed matrix between 1-forms and 2-forms $\mathbf{V}^1$ is given as $$\mathbf{V}^1_{ij,abc} = \langle \nabla \times \psi_{ij}, \psi_{abc} \rangle = 2 \sum_{(abc)} \tr \left( \tr \left( \mathbf{I}^{1,4}_{a,ibjc} - \mathbf{I}^{1,4}_{a,icjb} \right) \right) .$$
% \end{corollary}
\begin{corollary}
When the function $\psi_0=~1$ is added to the GRBF zero-forms,
% in our competed \textcolor{red}{space of 1-forms $\widehat{V}^1$}, 
the mass matrix of 1-forms $\widehat{~{M}}^1$ becomes
$$ \widehat{~{M}}^1 = \begin{bmatrix} ~S^0 & ~D^0 \\ {~D^0}^\intercal & ~M^1 \end{bmatrix},$$
and the mixed matrix $\widehat{~{D}}^0$ becomes
$$ \widehat{~{D}}^0 = \begin{bmatrix} ~{S}^0 & {~{D}^0} \end{bmatrix}.$$
Similar expressions can be constructed for $\widehat{~{M}}^2,\,\widehat{~{D}}^1$, etc.
\end{corollary}
% }\color{black}

Using these corollaries, we can formulate comparable machine-learning problems as in the first three examples in this section; the effectiveness of our mixed GRBF scheme is tested on the manufactured problem 
\begin{equation*} \begin{split}
~F - \nabla u &=  ~G \qquad \text{in } \Omega = \R^3,\\
\nabla \cdot ~F &= f \qquad \text{ in } \Omega, \\
\end{split} \end{equation*}
where $~G(~x) = \exp\left( -\nn{~x - \frac{1}{2} ~{1}}^2\right) ~{1}$ and $f(~x) = \exp\lr{-\nn{~x - \frac{1}{2}~{1}}^2} \left(3 - 2 ~x^\intercal~{1}\right)$.  Note that this problem has an analytic solution of 
$u(~x) = 0$ and $~F(~x) = ~G(~x).$
%We note that while the true solution $u = 0$, 
On the other hand, the exact solution in our approximation space of GRBFs is inherently nonzero due to the projection error of $f$ into the space $V_0 = \mathrm{span}\{\psi_i\}$; this particular problem is chosen to highlight approximation capabilities for problems where $~F \notin \text{im}(\nabla)$, which drives Darcy problems and other mixed-form discretizations of Poisson type.
% such mixed-form discretizations of these types of Darcy problems. 
Moreover, with this choice of $~G$, we know that there exist Gaussian RBF shape parameters so that the projection error of $~G$ into the space $V_1 = \mathrm{span}\{\psi_{ij}\}$ \emph{is} exactly 0, although this is almost surely not the case at initialization. This is due to the fact that, for Gaussians $\varphi_0(~x) = \exp(-\frac{1}{2}\nn{~x}^2)$ and $\varphi_1(~x) = \exp( - \frac{1}{2}\nn{~x-~{1}}^2)$, $~G$ is expressible as
\begin{equation*} \begin{split}
~G(~x) &= \exp\left( -\nn{~x - \frac{1}{2} ~{1}}^2\right) ~{1}\\
&= \exp\left( -\nn{~x - \frac{1}{2} ~{1}}^2\right) \left( \left( ~{1} - ~x \right) + ~x \right) \\
&= \exp\left(-\frac{1}{4} ~{1}^\intercal~{1}\right) \exp\left( - \frac{1}{2}\nn{~x}^2\right) \exp\left( - \frac{1}{2}\nn{~x - ~{1}}^2\right) \left( ~{1} - ~x \right) \\
&\qquad -  \exp\left(-\frac{1}{4} ~1^\intercal~1\right) \exp\left( - \frac{1}{2}\nn{~x}^2\right) \exp\left( - \frac{1}{2}\nn{~x - ~1}^2\right) \left( - ~x \right) \\
&=\exp\left(-\frac{d}{4} \right) \left( \varphi_0(~x) \nabla \varphi_1(~x) - \varphi_1(~x) \nabla \varphi_0(~x) \right).
\end{split} \end{equation*}

Table \ref{table:results-whitney} shows the results of training with a small number of RBFs. Due to the approximation results in Section~\ref{sec:approx}, we expect to see improvement in $u_h$ both with basis calibration (i.e., training) and with increasing GRBF basis size. 
% both with and without training (i.e. just from solving in the initial basis). 
Conversely, since there are no analogous guarantees for the one-form approximation space containing $~F_h$, we expect $~F_h$ to improve only with training.
% for $~F$ we expect for error to improve with training - the RBF approximation results provided earlier hold for approximation in $~V^0$ but not necessarily in $~V^1$, and our 1-form basis functions are not strictly RBFs.
% The expected improvement is generally observed, although the approximation error of $u$ with $n=8$ outperforms expectations on this example with the provided initializations.

\begin{table}[htbp!]
\centering
\begin{tabular}{r|ccc|ccc}
\hline
& \multicolumn{3}{c}{Solve} & \multicolumn{3}{c}{Train}\\

$N$ & Total Error & MSE $u$ & MSE $~F$ & Total Error & MSE $u$ & MSE $~F$ \\
\hline
8  & $9.4753 \times 10^{-4}$ & $4.1185 \times 10^{-4}$ & $5.3368 \times 10^{-4}$
   & $1.7717 \times 10^{-4}$ & $1.4227 \times 10^{-7}$ & $1.7703 \times 10^{-4}$ \\
16 & $9.3181 \times 10^{-4}$ & $3.9074 \times 10^{-4}$ & $5.4107 \times 10^{-4}$ 
   & $6.3768 \times 10^{-4}$ & $1.5383 \times 10^{-4}$ & $6.2229 \times 10^{-4}$ \\
32 & $1.0108 \times 10^{-3}$ & $4.4464 \times 10^{-4}$ & $5.6617 \times 10^{-4}$ 
   & $9.9938 \times 10^{-5}$ & $1.2075 \times 10^{-5}$ & $8.7863 \times 10^{-5}$\\
\hline
\end{tabular}
\caption{Relative data fit MSE for Problem 4, for an increasing number of Gaussian RBFs. \label{table:results-whitney}}
\end{table}

\section{Discussion and Conclusions}
\label{sec:conclusion}

Of interest is using this scheme to build surrogate models for quantities governed by Boltzmann-type equations, as these problems are posed in high-dimensions on unbounded domains. This focus remains an active area of ongoing work which we are currently investigating.

We remark that as $h \rightarrow 0$, this scheme effectively becomes a particle-based scheme using standard RBF collocation methods, and provides a natural framework to couple particle methods such as RKPM \cite{chen1996reproducing, behzadan2011unified} with more traditional variational formualations. Compared to particle-based collocation methods, due to the treatment of the boundary terms in our bilinear form, our scheme improves over e.g. Kansa's method \cite{kansa1990multiquadrics} since we guarantee that the resulting solution solves the underlying governing equations even while using a Gaussian kernel, and similarly we have the ability to construct stability estimates, which are not always provided with RBF interpolation methods.

\section*{Acknowledgments}
Sandia National Laboratories is a multimission laboratory managed and operated by National Technology \& Engineering Solutions of Sandia, LLC, a wholly owned subsidiary of Honeywell International Inc., for the U.S. Department of Energy’s National Nuclear Security Administration under contract DE-NA0003525. This paper describes objective technical results and analysis. Any subjective views or opinions that might be expressed in the paper do not necessarily represent the views of the U.S. Department of Energy or the United States Government. This article has been co-authored by an employee of National Technology \& Engineering Solutions of Sandia, LLC under Contract No. DE-NA0003525 with the U.S. Department of Energy (DOE). The employee owns all right, title and interest in and to the article and is solely responsible for its contents. The United States Government retains and the publisher, by accepting the article for publication, acknowledges that the United States Government retains a non-exclusive, paid-up, irrevocable, world-wide license to publish or reproduce the published form of this article or allow others to do so, for United States Government purposes. The DOE will provide public access to these results of federally sponsored research in accordance with the DOE Public Access Plan https://www.energy.gov/downloads/doe-public-access-plan.  The work of all authors is  supported by the U.S. Department of Energy, Office of Advanced Computing Research under the ``Scalable and Efficient Algorithms - Causal Reasoning, Operators, Graphs and Spikes" (SEA-CROGS) project.  The authors would like to acknowledge Chris Eldred (Sandia National Laboratories) and Xiaozhe Hu (Tufts University) for their insightful conversations regarding exterior calculus, Owen Davis (Sandia National Laboratories) for highlighting the connection to RBFNets, and Jason Torchinsky (Sandia National Laboratories) for fruitful discussions on training approaches for adaptive RBFs.

\bibliographystyle{siamplain}
\bibliography{references}

\section{Appendix A: Facts and Omitted Proofs}\label{app:proofs}

Here we record the statements and proofs that were omitted from the main body.

\subsection{Smooth Extensions to Unbounded Domains}
The below theorem, as stated in \cite{brenner2008mathematical}, provides for building a smooth extension of a function in a Sobolev space on a bounded domain to a comparable function on an unbounded domain. For a complete proof, see \cite{stein1970singular}.
\begin{theorem}
\label{thm:extension}
Let $f \in W^k_p(\Omega)$ for $\Omega \subset \R^n$ an open bounded domain with Lipschitz boundary. Then, there exists an extension mapping $E: W^k_p(\Omega) \rightarrow W^k_p(\R^n)$ such that $E f \lvert_\Omega = f$ for any $f \in W^k_p(\Omega)$ and such that $$\norm{E f}_{W^k_p(\R^n)} \leq c \norm{f}_{W^k_p(\Omega)},$$ where $c$ is independent of $f$.
\end{theorem}
To use this result, we modify this extension to create second extension that decays to zero away from $\Omega$ while maintaining smoothness. Let $\sigma_\Omega$ be an infinitely smooth cutoff function, given by $\sigma_\Omega = \mathbf{1}_{\Omega} * \phi_\epsilon$, such that $\sigma_\Omega(x) = 1$ for $x \in \Omega$, and that $\sigma_\Omega(x) = 0$ for $x \in \R^n$ where $\inf_{y \in \Omega} \norm{x-y}_2 < \epsilon$. Thus, we can extend any function $f$ specified on $\Omega$ to all of $\R^n$ in such a way so that the function $g \in W^k_p(\R^n)$ given by $g(x) := \sigma_\Omega(x) Ef(x)$ agrees with $f$ on $\Omega$, decays to zero outside of $\Omega$, and is sufficiently smooth. Moreover, by such a construction, we can make the integral of the smooth extension \textit{outside} of $\Omega$ to be negligibly small by controlling the width of the mollifier $\epsilon$, i.e. $$\left \lvert \int_{\R^n \backslash \Omega } f \right \rvert < \epsilon.$$

\subsection{Product of Gaussian PDFs}
\begin{proof}[Proof of Lemma~\ref{lem:gaussprod}]
This is a moderately involved induction on formulae for product of two Gaussians.  First, recall the Searle and Sherman-Morrison-Woodbury matrix inversion formulae (see, e.g., \cite{petersen2008matrix}),
\begin{align*}
    \lr{~A+~B}^{-1} &= ~A^{-1}\lr{~A^{-1}+~B^{-1}}^{-1}~B^{-1}, \\
    \lr{~A+~B}^{-1} &= ~A^{-1} - ~A^{-1}\lr{~A^{-1}+~B^{-1}}^{-1}~A^{-1}.
\end{align*}
The base case of the induction ($k=2$) proceeds as follows.  The product of two Gaussian probability densities (c.f. \cite{petersen2008matrix}) yields the normalization constant, mean, and covariance given by
\begin{align*}
    z_{:2} &= \lr{2\pi}^{-d/2}\det\lr{~C_1+~C_2}^{-1/2}\exp\lr{-\frac{1}{2}\lr{~m_1-~m_2}^\intercal\lr{~C_1+~C_2}^{-1}\lr{~m_1-~m_2}}, \\
    ~C_{:2}^{-1} &= ~C_1^{-1}+~C_2^{-1}, \\
    ~C_{:2}^{-1}~m_{:2} &= ~C_1^{-1}~m_1+~C_2^{-1}~m_2.
\end{align*}
Applying the Searle formula, the determinant term can be rewritten as
\[\det\lr{~C_1+~C_2}^{-1/2} = \lr{\det{\lr{~C_1+~C_2}}^{-1}}^{1/2} = \det\lr{~C_1^{-1}~C_{:2}~C_2^{-1}}^{-1/2} = \left(\frac{\det~C_{:2}}{\det~C_1\det~C_2}\right)^\frac{1}{2}.\]
Similarly, applying both inversion formulae to the exponential term yields
\begin{align*}
    \lr{~m_1-~m_2}^\intercal&\lr{~C_1+~C_2}^{-1}\lr{~m_1-~m_2} \\ 
    &= ~m_1^\intercal\lr{~C_1+~C_2}^{-1}~m_1 + ~m_2^\intercal\lr{~C_1+~C_2}^{-1}~m_2 - 2\,~m_1^\intercal\lr{~C_1+~C_2}^{-1}~m_2 \\
    &= ~m_1^\intercal\lr{~C_1^{-1}-~C_1^{-1}~C_{:2}~C_1^{-1}}~m_1 + ~m_2^\intercal\lr{~C_2^{-1}-~C_2^{-1}~C_{:2}~C_2^{-1}}~m_2 - 2\,~m_1^\intercal\lr{~C_1^{-1}~C_{:2}~C_2^{-1}}~m_2 \\
    &= ~m_1^\intercal~C_1^{-1}~m_1 + ~m_2^\intercal~C_2^{-1}~m_2 - ~m_1^\intercal~C_1^{-1}~C_{:2}\lr{~C_{:2}^{-1}~m_{:2}-~C_2^{-1}~m_2}\\
    &\qquad- ~m_2^\intercal~C_2^{-1}~C_{:2}\lr{~C_{:2}^{-1}~m_{:2}-~C_1^{-1}~m_1} - 2\,~m_1^\intercal\lr{~C_1^{-1}~C_{:2}~C_2^{-1}}~m_2 \\
    &= ~m_1^\intercal~C_1^{-1}~m_1 + ~m_2^\intercal~C_2^{-1}~m_2 - ~m_{:2}^\intercal~C_{:2}~m_{:2},
\end{align*}
establishing the Lemma when $k=2$.  Supposing it continues to hold for $k=n-1$, it follows that the product distribution becomes
\[\prod_{i=1}^n \phi_i = z_{:n-1}\phi_{:n-1}\phi_n = z_{:n}\phi_{:n}\]
where $\phi_{:n}$ is a density corresponding to the distribution $\mathcal{N}\lr{~m_{:n},~C_{:n}}$ and $z_{:n-1},~C_{:n-1},~m_{:n-1}$ are determined.  Applying the base case $k=2$ to $z_{:n-1}\phi_{:n-1}\phi_n$ then yields 
\begin{align*}
    z_{:n} &= z_{:n-1} \lr{2\pi}^{-\frac{d}{2}}\lr{\frac{\det~C_{:n}}{\det~C_{:n-1}\det~C_n}}^{\frac{1}{2}}\exp\lr{-\frac{1}{2}\lr{~m_{:n-1}^\intercal~C_{:n-1}^{-1}~m_{:n-1}+~m_{n}~C_n^{-1}~m_n-~m_{:n}~C_{:n}^{-1}~m_{:n}}} \\
    &= \lr{2\pi}^{-\frac{(n-2)d}{2}} \lr{\frac{\det~C_{:n-1}}{\prod_{i=1}^{n-1} \det~C_i}}^{\frac{1}{2}}\exp\lr{-\frac{1}{2}\lr{\sum_{i=1}^{n-1} ~m_i^\intercal~C_i^{-1}~m_i - ~m_{:n-1}^\intercal~C_{:n-1}^{-1}~m_{:n-1}}} \\
    &\qquad \cdot \lr{2\pi}^{-\frac{d}{2}}\lr{\frac{\det~C_{:n}}{\det~C_{:n-1}\det~C_n}}^{\frac{1}{2}}\exp\lr{-\frac{1}{2}\lr{~m_{:n-1}^\intercal~C_{:n-1}^{-1}~m_{:n-1}+~m_{n}~C_n^{-1}~m_n-~m_{:n}~C_{:n}^{-1}~m_{:n}}} \\
    % &= \lr{2\pi}^{-\frac{(n-2)d+d}{2}}\lr{\frac{\det~C_{:n-1}\det~C_{:n}}{\det~C_{:n-1}\prod_{i=1}^{n-1}\det~C_i\det~C_n}}^{\frac{1}{2}}\exp\lr{-\frac{1}{2}\lr{\sum_{i=1}^n ~m_i^\intercal~C_i^{-1}~m_i - ~m_{:n}^\intercal~C_{:n}^{-1}~m_{:n}}} \\
    &= \lr{2\pi}^{-\frac{(n-1)d}{2}} \lr{\frac{\det~C_{:n}}{\prod_{i=1}^n \det~C_i}}^{\frac{1}{2}}\exp\lr{-\frac{1}{2}\lr{\sum_{i=1}^n ~m_i^\intercal~C_i^{-1}~m_i - ~m_{:n}^\intercal~C_{:n}^{-1}~m_{:n}}},
\end{align*}
where the inductive hypothesis was used to expand $z_{:n-1}$ in the second line. This completes the inductive step and establishes the Lemma in the general case.
\end{proof}

\subsection{Polynomial Moments under Gaussian Measures}
% For a random vector $\mathbf{x} \in \mathbb{R}^n$, denote $\mathbf{x}^{\otimes d} := \mathbf{x} \otimes \mathbf{x} \otimes \cdots \otimes \mathbf{x}$ be the outer (tensor) product of $\mathbf{x}$ against itself $d$ times.
From \cite{pereira2022tensor}, we have the following lemma.
\begin{lemma}\label{lem:gaussian_moment} (From \cite{pereira2022tensor})
Let $~X \sim \mathcal{N}(~{m}, ~{C})$ be a random variable with probability density $\phi(~x)$. Then,
\begin{equation*}
\mathbb{E}( ~X^{\otimes d} ) = \int_{\mathbb{R}^n} ~x^{\otimes d} \, \phi(~x)\, d~x= \sum_{k=0}^{\left \lfloor d/2 \right \rfloor} C_{d,k}\,\mathrm{sym}\left( ~{m}^{\otimes(d-2k)} \otimes ~C^{\otimes k} \right), 
\end{equation*}
where, for any tensor $~T\in\mathbb{R}^{\otimes d}$ and denoting the permutation group on $d$ elements by $P_d$,
\begin{equation*}
    C_{d,k} = \begin{pmatrix} d \\ 2k \end{pmatrix} \frac{ (2k)! }{k! 2^k}, \quad \mathrm{sym}(~T) = \frac{1}{d!}\sum_{\sigma\in P_d} T_{\sigma_1...\sigma_d}.
\end{equation*}
\end{lemma}
Note that when $d=2$, this reduces to \Cref{lem:cov}.
% In the spirit of \cite{bader2006algorithm}, define the \textit{even-modal tensor product} $\times_E$ as follows: for tensor $A$ of size $n_1 \times n_2 \times \cdots \times n_{2p}$ and $B$ of size $ n_2 \times n_4 \times n_6 \times \cdots \times n_{2p}$, the even-mode tensor product $A \times_E B$ is given as $$ (A \times_E B)_{i_1,\dots,i_{2p-1}} = \sum_{i_2} \sum_{i_4} \cdots \sum_{i_{2p}} A_{i_1,i_2,\dots,i_{2p}} B_{i_2,i_4,i_6\dots,i_{2p}}. $$
Using Lemma \ref{lem:gaussian_moment}, we have the following result:

\begin{proof}[Proof of Theorem~\ref{thm:quadrature}]
The proof is separated into two parts.  First, we prove the result in terms of polynomial moments.  By the definition of the gradient of Gaussian functions, the integrand can be expressed as
\begin{equation*} \begin{split}\left( \prod_{a=1}^\alpha \phi_{i_a}(~{x}) \right) \left( \bigotimes_{b=1}^\beta \nabla \phi_{j_b}(~{x}) \right)
&= \left( \prod_{a=1}^\alpha \phi_{i_a}(~{x}) \right) \left( \prod_{b=1}^\beta \phi_{j_b}(~{x}) \right)  \left( \bigotimes_{b=1}^\beta ~{C}_{j_b}^{-1}\left( ~{m}_{j_b} - ~{x} \right) \right) \\
&= Z \Phi(~{x})  \left( \bigotimes_{b=1}^\beta ~{C}_{j_b}^{-1}\left( ~{m}_{j_b} - ~{x} \right) \right) \\
\end{split} \end{equation*}
where $\Phi$ is the probability distribution function for $~{X}$.
Integrating the above expression and making use of the fact that integrals against $\Phi$ are expectations, it follows that
\begin{equation*} \begin{split}
~{I}^{\alpha,\beta}_{i_1,\dots,i_\alpha;j_1,\dots,j_\beta} 
&= Z \, \mathbb{E} \left[ \bigotimes_{b=1}^\beta ~{C}_{j_b}^{-1}(~{m}_{j_b} - ~{X} ) \right] = Z \, \left( \bigotimes_{b=1}^\beta ~{C}_{j_b}^{-1} \right) \times_E \mathbb{E}\left[ \bigotimes_{b=1}^\beta ( ~{m}_{j_b} - ~{X}) \right],
\end{split} \end{equation*}
where the expectation is taken with respect to the distribution of $~X$.
We now expand the outer product inside the expectation into each of its constituent terms; each of these terms has a number of modes resulting from outer products of $~{X}$ with itself, and the remaining modes consist of means $~{m}_{j_b}$. In this expansion, we can group terms by the number of means $~m_{j_b}$ that appear in each expression.  It follows that
% $\mathbf{X}$'s that appear in each expression. Denoting the number of means $\mathbf{m}_{j_b}$ that appear with $k$, and thus the number of appearances of $\mathbf{X}$'s (i.e. the moment of $\mathbf{X}$ that needs to be computed in each term) with $\beta-k$, it follows that
\begin{equation*}
\begin{split}
    \mathbb{E}\left[ \bigotimes_{b-1}^\beta (~{m}_{j_b} - ~{X}) \right] 
    % &= \mathbb{E} \left[ \sum_{k=0}^\beta \sum_{\substack{J \in 2^{\lvert \beta \rvert} \\ \lvert J \rvert = k}} (-1)^{\beta-k} \left( \left( \bigotimes_{b=1}^k \mathbf{m}_{j_b} \right) \otimes \mathbf{X}^{\otimes (\beta-k)} \right)_{\sigma_{\vec{j}}} \right] \\
    &= \mathbb{E} \left[ \sum_{J\in 2^{\lvert \beta \rvert}} (-1)^{\beta-|J|} \left( \left( \bigotimes_{b=1}^{|J|} ~{m}_{j_b} \right) \otimes ~{X}^{\otimes (\beta-|J|)} \right)_{\sigma_{J}} \right], \\
    &= \sum_{J\in 2^{\lvert \beta \rvert}} (-1)^{\beta-|J|} \left( \left( \bigotimes_{b=1}^{|J|} ~{m}_{j_b} \right) \otimes ~M^{\beta-|J|} \right)_{\sigma_{J}},
\end{split}
\end{equation*}
where $~M^{\beta-|J|} = \mathbb{E}\left[~{X}^{\otimes (\beta-|J|)}\right]$ and the last equality follows by linearity of the expectation.
% Now, by linearity of expectation,
% \begin{equation*} \begin{split}
% \mathbb{E}\left[ \bigotimes_{b=1}^\beta (\mathbf{m}_{j_b} - \mathbf{X}) \right] 
% &= \sum_{k=0}^\beta \sum_{\substack{\vec{j} \in 2^{\lvert \beta \rvert} \\ \lvert \vec{j} \rvert = k}} (-1)^{\beta-k} \left( \left( \bigotimes_{b=1}^k \mathbf{m}_{j_b} \right) \otimes \mathbb{E} \left[  \mathbf{X}^{\otimes (\beta-k)} \right] \right)_{\sigma_{\vec{j}}} \\
% &= \sum_{k=0}^\beta \sum_{\substack{\vec{j} \in 2^{\lvert \beta \rvert}\\ \lvert \vec{j} \rvert = k}} (-1)^{\beta-k} \left( \left( \bigotimes_{b=1}^k \mathbf{m}_{j_b} \right) \otimes \mathbf{M}^{\beta-k}  \right)_{\sigma_{\vec{j}}}.
% \end{split} \end{equation*}
% where $2^{\lvert \beta \rvert}$ is the power set of the set $\{1,\dots,\beta \}$; $\vec{j}$ is a multi-index representing a subset of the power set $2^{\lvert \beta \rvert}$; and $\sigma_{\vec{j}}$ is the unique permutation of $\beta$ elements that maps the first $k = \lvert \vec{j} \rvert$ elements to positions $j_1,\dots,j_k$, and the rest of the elements \textit{in order} to positions $j_{k+1},\dots,j_\beta$.
Therefore,
\begin{equation*}
    ~{I}^{\alpha,\beta}_{i_1,\dots,i_\alpha; j_1,\dots,j_\beta} = Z \left( \bigotimes_{b=1}^\beta ~{C}_{j_b}^{-1} \right) \times_E  \sum_{J\in 2^{\lvert \beta \rvert}} (-1)^{\beta-|J|} \left( \left( \bigotimes_{b=1}^{|J|} ~{m}_{j_b} \right) \otimes ~M^{\beta-|J|} \right)_{\sigma_{J}}.
\end{equation*}
% \left(  \sum_{k=0}^\beta \sum_{\substack{\vec{j} \in 2^{\lvert \beta \rvert} \\ \lvert \vec{j} \rvert = k}} (-1)^{\beta-k} \left( \left(  \bigotimes_{b=1}^k \mathbf{m}_{j_b} \right) \otimes \mathbf{M}^{\beta-k} \right)_{\sigma_{\vec{j}}} \right)
The linearity of the even mode tensor product completes the expression in terms of polynomial moments.

The next goal is to prove the expression in terms of Gauss-Hermite quadrature.  For ease of notation, define the tensor-valued polynomial 
\begin{equation*}
    ~{p}(~{x}) = \bigotimes_{b=1}^\beta ~{C}^{-1}_{j_{b}} \left(~{m}_{j_b} - ~{x} \right).
\end{equation*}
By the product Lemma~\ref{lem:gaussprod}, it follows that
\begin{equation*}
\begin{split}
    ~{I}^{\alpha,\beta}_{i_1,\dots,i_\alpha;j_1,\dots,j_\beta} &= \int_{\mathbb{R}^n} \left( \prod_{a=1}^\alpha \phi_{i_a}(~{x}) \right) \left( \bigotimes_{b=1}^\beta \nabla \phi_{j_b}(~{x}) \right) \, d~{x} \\
    &= \int_{\mathbb{R}^n} ~{p}(~{x}) \left( \prod_{a=1}^\alpha \phi_{i_a}(~{x}) \right) \left( \prod_{b=1}^\beta \phi_{j_b}(~{x}) \right) \, d~{x} = Z \int_{\mathbb{R}^n} ~{p}(~{x}) \phi(~{x}) \, d~{x},
\end{split} 
\end{equation*}
where $\phi=\phi_{i_1,...,i_\alpha,j_1,...,j_\beta}$ denotes the probability density function of the product with normalization factor $Z$.  Expanding out the definition for the probability density function $\phi$ and applying a change of variables $\widehat{~{x}} = \frac{1}{\sqrt{2}} ~{C}^{-\frac{1}{2}}\left(~{x}-~{m}\right)$ then yields
\begin{equation*}
\begin{split}
    ~{I}^{\alpha,\beta}_{i_1,\dots,i_\alpha;j_1,\dots,j_\beta}
    &= Z \left(2 \pi\right)^{-\frac{n}{2}} \det(~{C})^{-\frac{1}{2}} \int_{\mathbb{R}^n} \exp \left(-\frac{1}{2}\left(~{x}-~{m}\right)^\intercal ~{C}^{-1}\left(~{x}-~{m} \right) \right) ~{p}(~{x}) \, d~{x} \\
    &= Z \left(2 \pi\right)^{-\frac{n}{2}} \det(~{C})^{-\frac{1}{2}} \int_{\mathbb{R}^n} \exp\left(-\widehat{~{x}}^\intercal \widehat{~{x}} \right) ~{p}\left(\sqrt{2}~{C}^{\frac{1}{2}}\widehat{~{x}} + ~{m}\right) \, \left( 2^{\frac{d}{2}} \det(~{C})^{\frac{1}{2}}\right) \, d\widehat{~{x}} \\
    &= Z\, \pi^{-\frac{n}{2}} \int_{\mathbb{R}^n} \exp\left(-\widehat{~{x}}^\intercal \widehat{~{x}} \right) ~{p}\left(\sqrt{2}~{C}^{\frac{1}{2}}\widehat{~{x}} + ~{m}\right) \, d\widehat{~{x}}.
\end{split} \end{equation*}
Applying Gauss-Hermite quadrature then provides the final result, since the quadrature was chosen as to integrate $~{p}$ exactly under a squared-exponential weight. 
\end{proof}

\newpage
\section{Supplemental Material}
Here we record additional experimental details, as well as alternative, direct calculations of the mass and stiffness matrices necessary for implementing the examples in the body.

\subsection{Additional Figures and Details}
\begin{figure}[htbp!]
\includegraphics[height=2in]{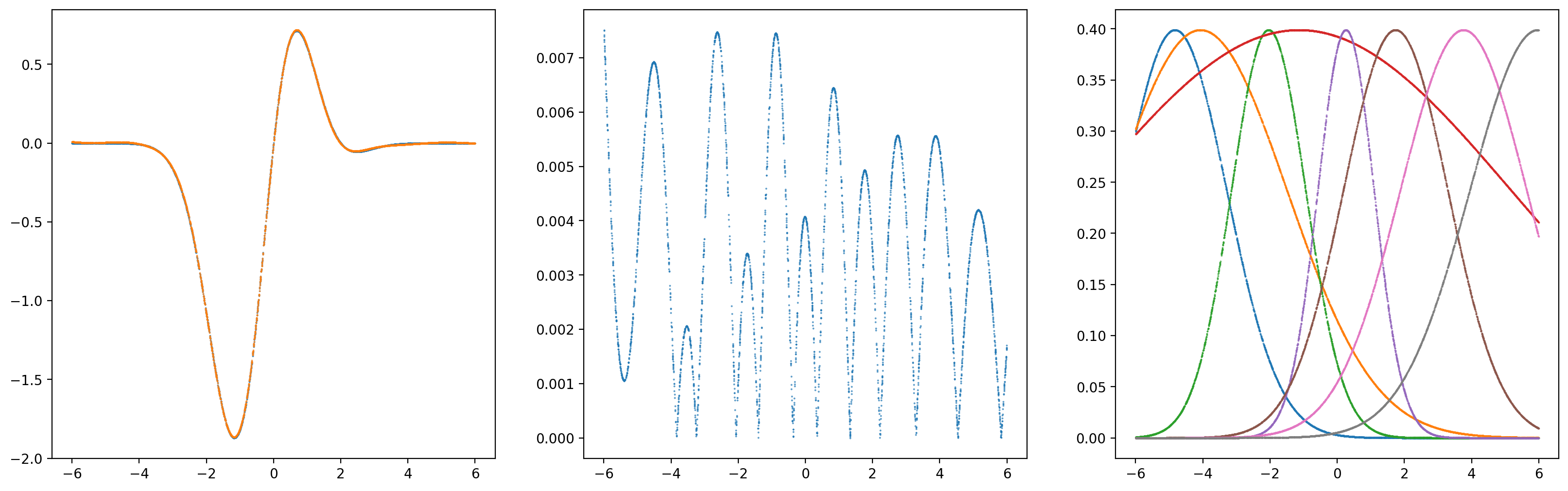} \\
\includegraphics[height=2in]{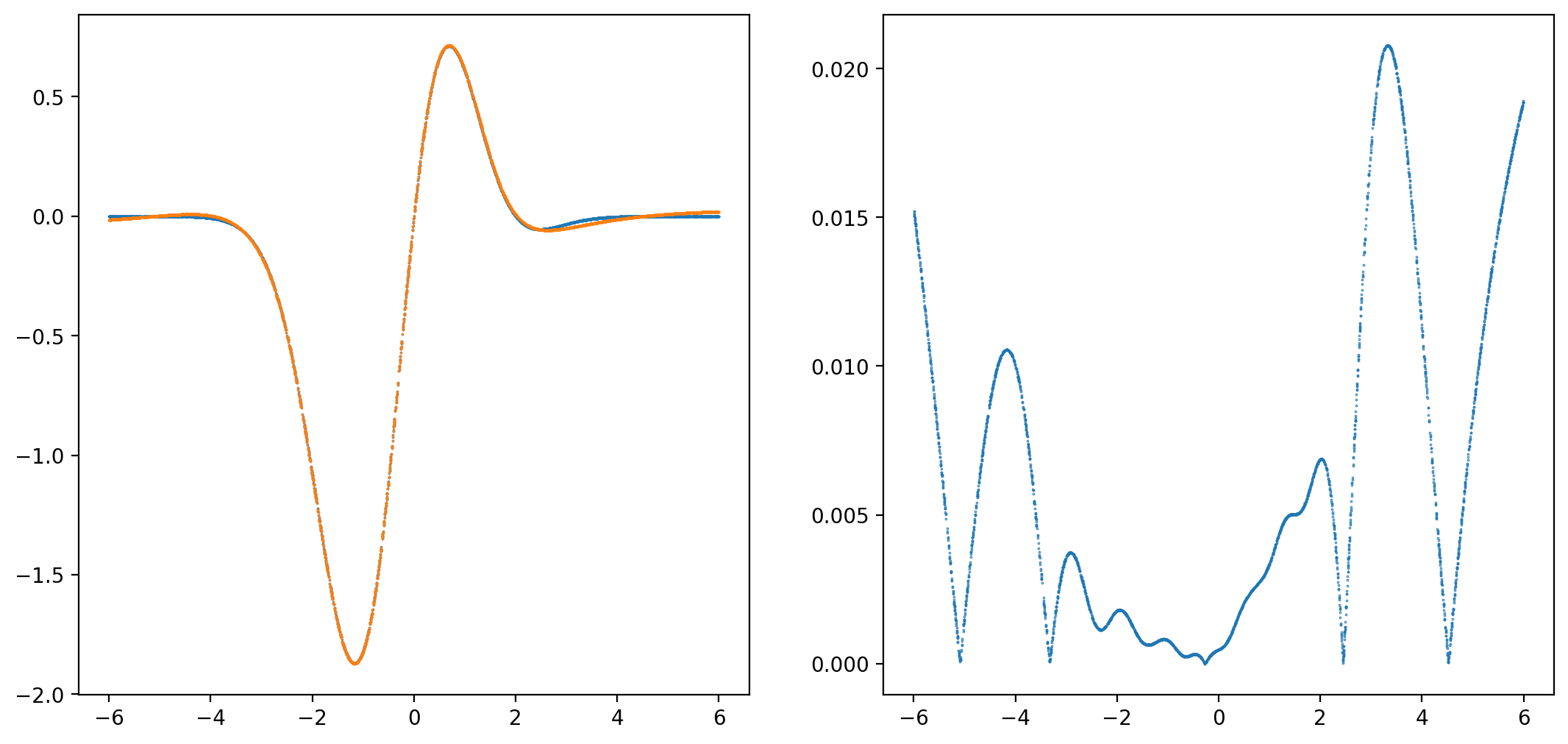} \\
\includegraphics[height=2in]{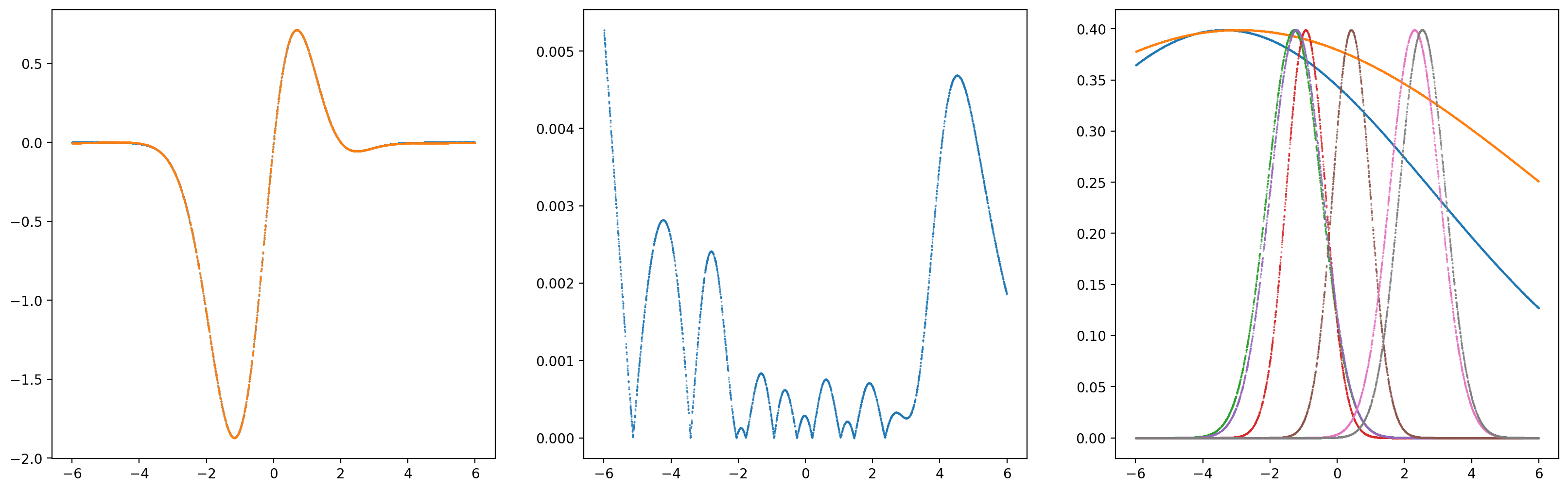} 
\caption{Solutions to Problem 1 for $N=8$ using our novel GRBF method (top), a standard PINN (middle), and an RBFNet PINN (bottom). From left to right: plot of true (blue) and predicted (orange) solutions; pointwise absolute error; scaled Gaussian basis found during training. \label{fig:prob1-train8-accuracy}}
\end{figure}

\begin{figure}[htbp!]
\includegraphics[height=2in]{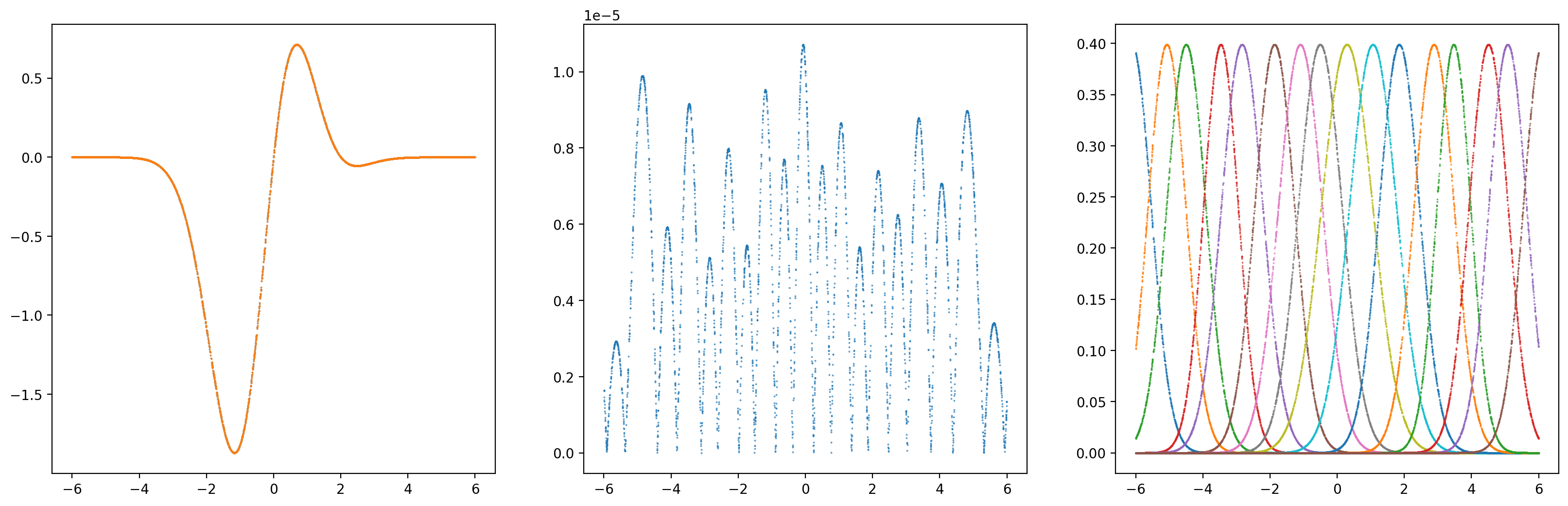} \\
\includegraphics[height=2in]{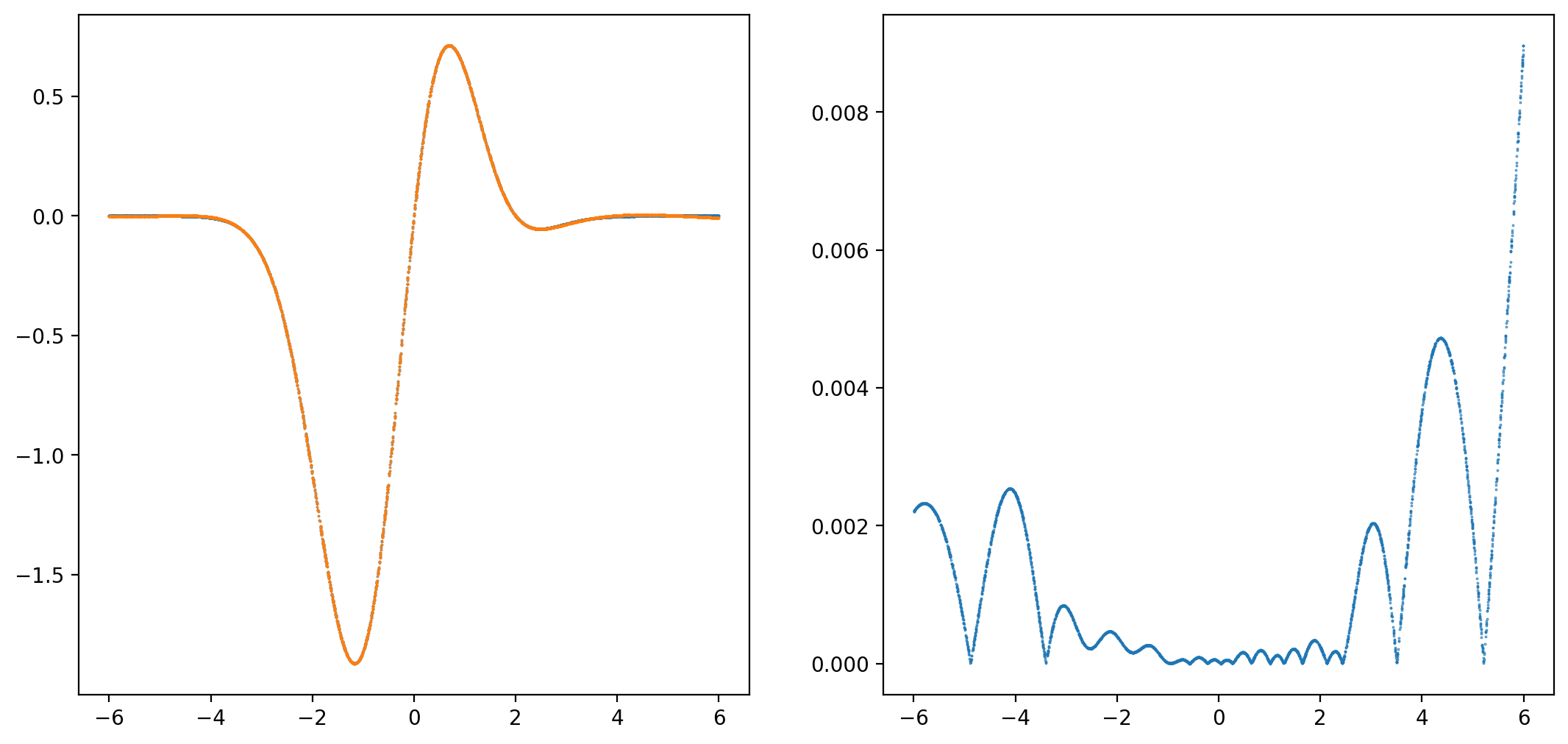} \\
\includegraphics[height=2in]{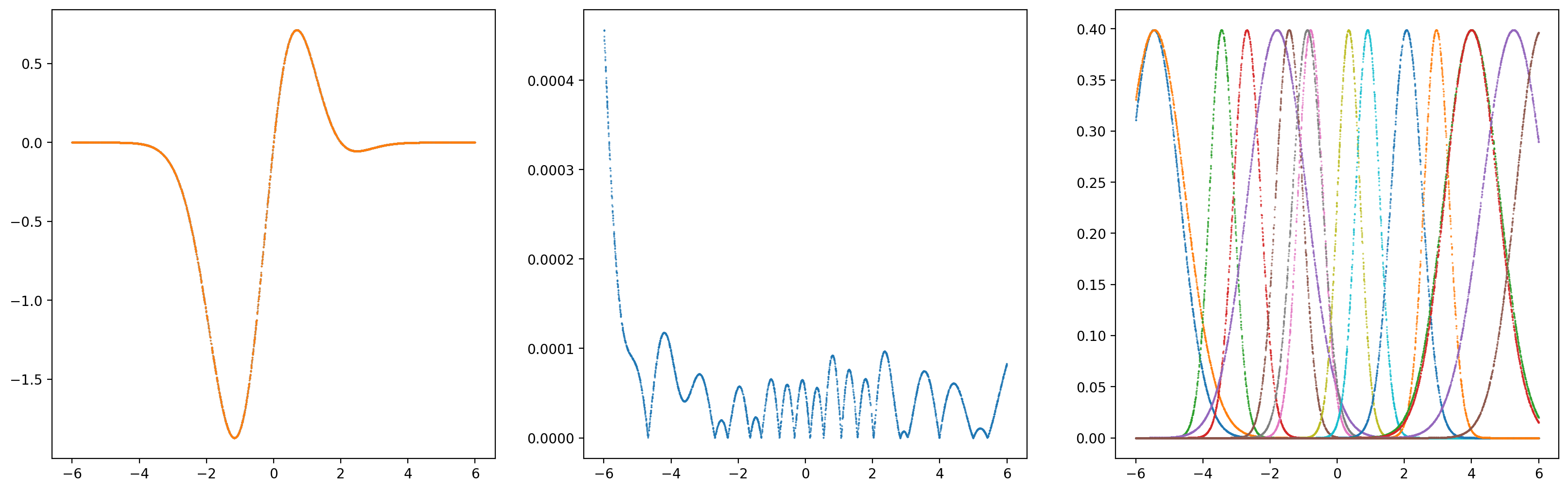} 
\caption{Solutions to Problem 1 for $N=16$ using our novel GRBF method (top), a standard PINN (middle), and an RBFNet PINN (bottom). From left to right: plot of true (blue) and predicted (orange) solutions; pointwise absolute error; scaled Gaussian basis found during training. \label{fig:prob1-train16-accuracy}}
\end{figure}

% \label{app:fig}
Figure \ref{fig:prob1-train8-accuracy} compares the solutions to Problem 1, and the recovered Gaussian bases, for $N=8$; Figure \ref{fig:prob1-train16-accuracy} does the same for $N=16$.

\subsection{Integrals of Trigonometric Functions against Gaussian Distributions}
% \label{app:trig}

We first present the 1D case here, of integrating trigonometric functions against a Gaussian probability distribution. For higher dimensions, similar formulae can be derived for both the product of univariate trigonometric functions by separability, $u$-substitution, and various trigonometric identities, and also for trig functions composed with affine functions by a similar approach.
\begin{lemma}
Let $X\sim \mathcal{N}(\mu, \sigma^2)$ be a normal random variable with PDF $\phi$. Then,
\begin{equation*} \begin{split}
\int_{-\infty}^{\infty} \sin (k x) \phi(x) \,dx 
&= \sin( k \mu ) \exp \left( -\frac{1}{2} k^2 \sigma^2 \right), \\ 
\text{and}\qquad \int_{-\infty}^{\infty} \cos (k x) \phi(x) \,dx &= \cos( k \mu ) \exp \left( -\frac{1}{2} k^2 \sigma^2 \right).
\end{split} \end{equation*}
\end{lemma}
\begin{proof}
We first make a few observations:
\begin{enumerate}
\item By symmetry and since $\sin$ is an odd function, we have $$\int_{-\infty}^\infty \sin(a x) \exp \left( - x^2 \right) \, dx = 0.$$
\item By the Cauchy Integral Theorem, 
\begin{equation*} \begin{split}
\int_{-\infty}^\infty \cos(a x) \exp \left( - x^2 \right) \, dx 
&= \Re \int_{-\infty}^{\infty} \exp\left( \mathrm{i} a x - x^2 \right) \, dx \\
&= \exp\left(-\frac{a^2}{4}\right) \Re \int_{-\infty}^\infty \exp \left( - \left( x - \frac{\mathrm{i}a}{2} \right)^2 \right) \, dx \\
&= \sqrt{\pi} \exp\left(-\frac{a^2}{4}\right).
\end{split} \end{equation*}
\item From fundamental trigonometric relationships,
\begin{equation*} \begin{split}
\sin (x + y) &= \sin(x) \cos(y) + \cos(x) \sin(y) \\
\cos (x + y) &= \cos(x) \cos(y) - \sin(x) \sin(y).
\end{split} \end{equation*}
\end{enumerate}
Putting these observations together,
\begin{equation*} \begin{split}
\int_{-\infty}^\infty \sin(k x) \phi(x) \, dx
&= \frac{1}{\sqrt{2 \pi} \sigma} \int_{-\infty}^\infty \sin( k x) \exp\left( - \frac{1}{2}\left( \frac{x - \mu}{\sigma} \right)^2 \right) \, dx \\
&= \frac{1}{\sqrt{\pi} } \int_{-\infty}^\infty \sin\left( \sqrt{2}\sigma k z + k \mu \right) \exp\left( - z^2 \right) \, dz \\
&= \frac{1}{\sqrt{\pi} } \int_{-\infty}^\infty \left( \sin\left( \sqrt{2}\sigma k z\right) \cos\left( k \mu \right) + \cos\left( \sqrt{2} \sigma k z \right) \sin\left(k \mu \right)\right) \exp\left( - z^2 \right) \, dz \\
&= \frac{1}{\sqrt{\pi} } \sin(k \mu) \int_{-\infty}^\infty \cos\left(\sqrt{2}\sigma k z \right) \exp(-z^2) \,dz \\
&= \sin(k \mu) \exp\left(-\frac{1}{2} \sigma^2 k^2 \right).
\end{split} \end{equation*}
Similarly,
\begin{equation*} \begin{split}
\int_{-\infty}^\infty \cos(k x) \phi(x) 
&= \frac{1}{\sqrt{2 \pi} \sigma} \int_{-\infty}^\infty \cos( k x ) \exp\left(-\frac{1}{2}\left(\frac{x - \mu}{\sigma}\right)^2 \right) \, dx \\
&= \frac{1}{\sqrt{\pi}} \int_{-\infty}^\infty \left( \cos\left( \sqrt{2}\sigma k z\right) \cos\left( k \mu \right) - \sin\left(\sqrt{2}\sigma k z \right) \sin\left( k \mu \right) \right) \exp\left( - z^2 \right) \, dz \\
&= \frac{1}{\sqrt{\pi}} \cos( k \mu ) \int_{-\infty}^{\infty} \cos \left( \sqrt{2}\sigma k z \right) \exp\left( - z^2 \right) \, dz \\
&= \cos( k \mu ) \exp\left(-\frac{1}{2} \sigma^2 k^2 \right).
\end{split} \end{equation*}
\end{proof}

\subsection{Closed Form Expressions for \texorpdfstring{$k$}-Form Moments}
% \label{app:handcalc}
Closed-form expressions for various moments of $k$-forms generated from a Gaussian basis are now developed. These results are generalized to $k$-forms of arbitrary order by Theorem~\ref{thm:quadrature} of the body.  For ease of notation, this Section drops the distinction between the random variable $~X$ and its values $~x$.

\begin{proof}[Direct Proof of Corollary~\ref{cor:exxt}]
    By linearity of the expectation, we have
\begin{align*}
    ~C &= \expec{(~x-~m)(~x-~m)^\intercal} = \expec{~x~x^\intercal}-~m\,\expec{~x^\intercal}-\expec{~x}~m^\intercal+~m~m^\intercal \\
    &= \expec{~x~x^\intercal} - 2~m~m^\intercal + ~m~m^\intercal = \expec{~x~x^\intercal}-~m~m^\intercal,
\end{align*}
as desired.
\end{proof}
\begin{proof}[Direct Proof of Corollary~\ref{cor:equadform}]
Since the mean and covariance of $~x$ are $~m, ~C$, it follows from Lemma~\ref{lem:cov} and linearity of the expectation that
\begin{align*}
    \expec{(~x-~b)^\intercal~A(~x-~a)} &= \IP{\expec{(~x-~b)(~x-~a)^\intercal}}{~A}{} = \IP{\expec{~x~x^\intercal-~b~x^\intercal-~x~a^\intercal+~b~a^\intercal}}{~A}{} \\
    &= \IP{~C+~m~m^\intercal-~b~m^\intercal-~m~a^\intercal+~b~a^\intercal}{~A}{} \\
    &= \IP{~C}{~A}{}+~m^\intercal~A~m+~b^\intercal~A~a-~b^\intercal~A~m-~m^\intercal~A~a,
\end{align*}
from which the result follows by grouping the last four terms.
\end{proof}

% \JONAS{
% Things to clean up:
% \begin{itemize}
% \item curly $~x$ vs $\mathbf{x}$ and similar typesetting (C's, M's, etc.)
% \item transpose $^\intercal$ vs $^T$ vs $^t$ (right now all three are used)
% \item do indices for 1-forms include $i<j$ or all $ij$?
% \item which formulations of 1-forms include the basis function $\mathbf{1}$?
% \end{itemize}
% }

\subsubsection{Preliminaries for one-forms}

We repeat the following lemmas stated previously in Section \ref{sec:method}.

\begin{lemma}\label{lem:cov}
    Let $~x\sim\mathcal{N}\lr{~m, ~C}$, then 
    \[\expec{~x~x^\intercal} = ~C+~m~m^\intercal.\]
\end{lemma}
\begin{proof}
    By linearity of the expectation, we have
\begin{align*}
    ~C &= \expec{(~x-~m)(~x-~m)^\intercal} = \expec{~x~x^\intercal}-~m\,\expec{~x^\intercal}-\expec{~x}~m^\intercal+~m~m^\intercal \\
    &= \expec{~x~x^\intercal} - 2~m~m^\intercal + ~m~m^\intercal = \expec{~x~x^\intercal}-~m~m^\intercal,
\end{align*}
as desired.
\end{proof}

\begin{lemma}\label{lem:quadform}
Let $~x\sim\mathcal{N}\lr{~m, ~C}$, $~a,~b$ be generic vectors, and $~A$ be a generic square matrix.  Then 
\[\expec{(~x-~b)^\intercal~A(~x-~a)} = \IP{~C}{~A}{} + (~m-~b)^\intercal~A(~m-~a).\]
\end{lemma}
\begin{proof}
Since the mean and covariance of $~x$ are $~m, ~C$, it follows from Lemma~\ref{lem:cov} and linearity of the expectation that
\begin{align*}
    \expec{(~x-~b)^\intercal~A(~x-~a)} &= \IP{\expec{(~x-~b)(~x-~a)^\intercal}}{~A}{} = \IP{\expec{~x~x^\intercal-~b~x^\intercal-~x~a^\intercal+~b~a^\intercal}}{~A}{} \\
    &= \IP{~C+~m~m^\intercal-~b~m^\intercal-~m~a^\intercal+~b~a^\intercal}{~A}{} \\
    &= \IP{~C}{~A}{}+~m^\intercal~A~m+~b^\intercal~A~a-~b^\intercal~A~m-~m^\intercal~A~a,
\end{align*}
from which the result follows by grouping the last four terms.
\end{proof}

\subsubsection{One-form mass matrix}
\label{subsec:one-forms}
Begin with $~x\in\mathbb{R}^d$ and Gaussian densities $\phi_i(~x)$, $1\leq i\leq n$, as zero-forms.  Then, it follows immediately that $\nabla\phi_i(~x) = -\phi_i(~x)~C^{-1}_i(~x-~m_i)$.  Defining $\phi_{0} = 1$ 
%\JONAS{[Note - need to be consistent with where we develop the function $\mathbf{1}$ in our basis; need to be consistent about whether we construct $i<j$ 1-forms or all $ij$ 1-forms]} 
and building Whitney one-forms using the set of functions $\{\phi_0,\phi_1,...,\phi_n\}$ leads to
(a proxy for) the basis Whitney 1-form,
\begin{align*}
    \psi_{ij} &= \phi_i\nabla\phi_j-\phi_j\nabla\phi_i = \phi_i(~x)\phi_j(~x)\lr{~C^{-1}_i(~x-~m_i)-~C^{-1}_j(~x-~m_j)}, \quad 0\leq i,j\leq n,
\end{align*}
where, obviously, $~m_0,~C_0^{-1}$ are identically zero. 
%\JONAS{What does it mean to take the inverse of a matrix that's identically 0?}  
The goal is to compute the mass matrix with $(n+1)^4$ entries
\begin{align*}
    M^{1}_{ij,ab} &= \int_\mathbb{R}^d \psi_{ij}\cdot\psi_{ab} \\
    &= \int_\mathbb{R}^d \phi_i(~x)\phi_j(~x)\phi_a(~x)\phi_b(~x)\lr{~C^{-1}_i(~x-~m_i)-~C^{-1}_j(~x-~m_j)}\cdot\lr{~C^{-1}_a(~x-~m_a)-~C^{-1}_b(~x-~m_b)} \\
    &= z_{ijab}\int_\mathbb{R}^d \phi_{ijab}(~x)\lr{~C^{-1}_i(~x-~m_i)-~C^{-1}_j(~x-~m_j)}\cdot\lr{~C^{-1}_a(~x-~m_a)-~C^{-1}_b(~x-~m_b)} \\
    &= z_{ijab}\,\mathbb{E}\left[\lr{~C^{-1}_i(~x-~m_i)-~C^{-1}_j(~x-~m_j)}\cdot\lr{~C^{-1}_a(~x-~m_a)-~C^{-1}_b(~x-~m_b)}\right]
\end{align*}
where $\mathbb{E}$ denotes expectation with respect to the Gaussian density $\phi_{ijab}$ and Lemma~\ref{lem:gaussprod} was used in the third equality.  Note that, letting $1\leq i,j,a,b\leq n$ temporarily, there are only 3 nontrivial cases here due to antisymmetry in $i,j$ and $a,b$ and pairwise symmetry in $ij,ab$. These can be represented by $\mathcal{M}_{ij,ab}, \mathcal{M}_{ij,a0}, \mathcal{M}_{i0,a0}$, and the same formulae can be used for all of them provided care is taken in the products of Gaussian densities.
% First, notice that $\phi_{ijab}(~x) := \phi_i(~x)\phi_j(~x)\phi_a(~x)\phi_b(~x)$ is a product of Gaussians, so $\phi_{ijab}\sim Z_{ijab}\,\mathcal{N}(~m_{ijab}, ~C_{ijab})$ where
% \begin{align*}
%     Z_{ijab} &= \lr{2\pi}^{-d/2}\det\lr{~C_{ija}+~C_b}^{-1/2}\exp\lr{-\frac{1}{2}\lr{~m_{ija}-~m_b}^\intercal\lr{~C_{ija}+~C_b}^{-1}\lr{~m_{ija}-~m_b}}, \\
%     ~C_{ijab} &= \lr{~C^{-1}_i+~C^{-1}_j+~C^{-1}_a+~C^{-1}_b}^{-1}\\
%     ~m_{ijab} &= ~C_{ijab}\lr{~C^{-1}_i~m_i+~C^{-1}_j~m_j+~C^{-1}_a~m_a+~C^{-1}_b~m_b}.
% \end{align*}
Using Lemma~\ref{lem:quadform}, 
% and letting $\mathbb{E}$ denote the expectation with respect to $\mathcal{N}\lr{~m_{ijab},~C_{ijab}}$ 
we have 
\begin{align*}
    M^1_{ij,ab} &= z_{ijab}\,\expec{(~x-~m_i)^\intercal~C_i^{-\intercal}~C^{-1}_a(~x-~m_a) + (~x-~m_j)^\intercal~C_j^{-\intercal}~C^{-1}_b(~x-~m_b)} \\
    &\quad- z_{ijab}\,\expec{(~x-~m_i)^\intercal~C_i^{-\intercal}~C^{-1}_b(~x-~m_b) + (~x-~m_j)^\intercal~C_j^{-\intercal}~C^{-1}_a(~x-~m_a)} \\
    &= z_{ijab}\IP{~C_i^{-\intercal}~C^{-1}_a+~C_j^{-\intercal}~C^{-1}_b-~C_i^{-\intercal}~C^{-1}_b-~C_j^{-\intercal}~C^{-1}_a}{~C_{ijab}}{}\\
    &\quad+ z_{ijab}\lr{\lr{~m_{ijab}-~m_i}^\intercal~C_i^{-\intercal}~C^{-1}_a\lr{~m_{ijab}-~m_a} + \lr{~m_{ijab}-~m_j}^\intercal~C_j^{-\intercal}~C^{-1}_b\lr{~m_{ijab}-~m_b}} \\
    &\quad- z_{ijab}\lr{\lr{~m_{ijab}-~m_i}^\intercal~C_i^{-\intercal}~C^{-1}_b\lr{~m_{ijab}-~m_b} + \lr{~m_{ijab}-~m_j}^\intercal~C_j^{-\intercal}~C^{-1}_a\lr{~m_{ijab}-~m_a}},
\end{align*}
in agreement with Corollary~\ref{cor:1formM}.  This can be written more compactly by defining the operators
\begin{align*}
\lr{~C^{-1}~m}^{T}_{ijab} = 
    \begin{pmatrix}
        ~C_i^{-1}\lr{~m_{ijab}-~m_i} \\ ~C_j^{-1}\lr{~m_{ijab}-~m_j}
    \end{pmatrix}, \quad
\lr{~C^{-1}~m}^{t}_{ijab} = 
    \begin{pmatrix}
        ~C_a^{-1}\lr{~m_{ijab}-~m_a} \\ ~C_b^{-1}\lr{~m_{ijab}-~m_b}
    \end{pmatrix}, \quad
~J = -
    \begin{pmatrix}
        ~0 & ~I \\ ~I & ~0
    \end{pmatrix}.
\end{align*}
Using this, it follows that 
\[M^1_{ij,ab} = z_{ijab}\IP{~C_{ijab}}{\lr{~C_i^{-1}-~C_j^{-1}}^\intercal\lr{~C_a^{-1}-~C_b^{-1}}}{} + z_{ijab}\,\lr{\lr{~C^{-1}~m}^{T}_{ijab}}^\intercal\lr{~I-~J}\lr{~C^{-1}~m}^{t}_{ijab}.\]

\subsubsection{Preliminaries for two-forms}
To compute the mass matrix of 2-forms by this technique, we need the following moderately involved Lemmata.

% \begin{remark}
%     While the Lemmas in this section are stated for symmetric matrices $~A,~B$, they are directly extensible to generic square matrices since $~x^\intercal~A~x = (1/2)~x^\intercal\lr{~A+~A^\intercal}~x$.
% \end{remark}

\begin{lemma}\label{lem:triform}
    Let $~x\sim\mathcal{N}\lr{~m, ~C}$, then for generic square matrix $~A$ it follows that
    \[\expec{\lr{~x^\intercal~A~x}~x} = \IP{~A}{~C}~m + ~C\lr{~A+~A^\intercal}~m + \lr{~m^\intercal~A~m}~m.\]
\end{lemma}
\begin{proof}
    Parameterize $~x = ~m + ~C^{\frac{1}{2}}~y$ for $~y\sim\mathcal{N}\lr{~0,~I}$.  Then denoting $~S = ~C^{\frac{1}{2}}$ for convenience, it follows for symmetric $~A$ that
    \begin{align*}
        \lr{~x^\intercal~A~x}~x &= \lr{~m^\intercal~A~m}~m + \lr{\lr{~m^\intercal~A~m}~S~y + 2\lr{~m^\intercal~A~S~y}~m} \\
        &\quad+ \lr{\lr{~y^\intercal~S~A~S~y}~m + 2~S~y~y^\intercal~S~A~m} + \lr{~y^\intercal~S~A~S~y}~S~y.
    \end{align*}
    Taking expectations, it follows that the terms of odd degree in $~y$ will drop out, leaving
    \begin{align*}
        \expec{\lr{~x^\intercal~A~x}~x} &= \lr{~m^\intercal~A~m}~m + \expec{~y^\intercal~S~A~S~y}~m + 2~S\,\expec{~y~y^\intercal}~S~A~m \\
        &= \lr{~m^\intercal~A~m}~m + 2\,~C~A~m + \IP{~A}{~C}~m,
    \end{align*}
since $~S^2=~C$ and the trace obeys a cyclic property.  Now, notice that for generic square $~A$ we have $~x^\intercal~A~x = ~x^\intercal~A^\intercal~x = \frac{1}{2}~x^\intercal\lr{~A+~A^\intercal}~x$, so that we conclude 
\begin{align*}
    \expec{\lr{~x^\intercal~A~x}~x} &= \frac{1}{2}\IP{~A+~A^\intercal}{~C} + ~C\lr{~A+~A^\intercal}~m + \frac{1}{2}\lr{~m^\intercal\lr{~A+~A^\intercal}~m}~m \\
    &= \IP{~A}{~C}~m + ~C\lr{~A+~A^\intercal}~m + \lr{~m^\intercal~A~m}~m,
\end{align*}
as desired.
\end{proof}

\begin{lemma}\label{lem:mean0quadprod}
    Let $~y\sim\mathcal{N}\lr{~0, ~I}$, then for generic square $~A,~B$ it follows that
    \[\expec{\lr{~y^\intercal~A~y}\lr{~y^\intercal~B~y}} = \IP{~A+~A^\intercal}{~B} + \lr{\tr\,~A}\lr{\tr\,~B}.\]
\end{lemma}
\begin{proof}
    Recall that the components $y_i$ of $~y$ are i.i.d. and $\expec{y_i}=\expec{y_i^3}=0$. It follows that only second order terms survive in the expectation, and since $\expec{y_i^4}=3$ we have, for symmetric $~A,~B$,
    \begin{align*}
        \expec{\lr{~y^\intercal~A~y}\lr{~y^\intercal~B~y}} &= \sum_{i=j} A_{ii}B_{ii}\expec{y_i^4} + \sum_{i\neq j}A_{ii}B_{jj}\expec{y_i^2y_j^2} + 2\sum_{i\neq j}A_{ij}B_{ij}\expec{y_i^2y_j^2} \\
        &= 3\sum_{i}A_{ii}B_{ii} + \sum_{i\neq j} A_{ii}B_{jj} + \sum_{i\neq j} A_{ij}B_{ij} = \sum_{i,j}A_{ii}B_{jj} + 2\lr{\sum_{i\neq j}A_{ij}B_{ij} + \sum_{i=j}A_{ij}B_{ij}} \\
        &= \sum_iA_{ii}\sum_jB_{jj} + 2\sum_{i,j}A_{ij}B_{ij} = \lr{\tr\,~A}\lr{\tr\,~B} + 2\IP{~A}{~B}.
    \end{align*}
    The result for generic $~A,~B$ now follows by symmetrizing both sides of the above.  More precisely, notice that $\tr\,~A = \tr\,~A^\intercal = \frac{1}{2}\,\tr\lr{~A+~A^\intercal}$, and similarly
    \begin{align*}
        \frac{1}{4}\IP{~A+~A^\intercal}{~B+~B^\intercal} &= \frac{1}{4}\lr{\IP{~A}{~B} + \IP{~A^\intercal}{~B^\intercal}} + \frac{1}{4}\lr{\IP{~A^\intercal}{~B} + \IP{~A}{~B^\intercal}} \\
        &= \frac{1}{2}\IP{~A+~A^\intercal}{~B},
    \end{align*}
    establishing the conclusion.  Note that this expression is also equivalent to $\frac{1}{2}\IP{~A}{~B+~B^\intercal}$.
\end{proof}

\begin{lemma}\label{lem:quadprod}
    Let $~x\sim\mathcal{N}\lr{~m, ~C}$, then for generic square $~A,~B$ it follows that
    \begin{align*}
        \expec{\lr{~x^\intercal~A~x}\lr{~x^\intercal~B~x}} &= \IP{\lr{~A+~A^\intercal}~C}{~C~B} + ~m^\intercal\lr{~A+~A^\intercal}~C\lr{~B+~B^\intercal}~m \\
        &\quad + \lr{~m^\intercal~A~m + \IP{~A}{~C}}\lr{~m^\intercal~B~m + \IP{~B}{~C}}.
    \end{align*}
\end{lemma}
\begin{proof}
    Parameterize $~x = ~m + ~C^{\frac{1}{2}}~y$ for $~y\sim\mathcal{N}\lr{~0,~I}$.  Then denoting $~S = ~C^{\frac{1}{2}}$ for convenience, it follows for symmetric $~A,~B$ that (grouping terms suggestively)
    \begin{align*}
        (~x^\intercal&~A~x) \lr{~x^\intercal~B~x} = \lr{~m+~S~y}^\intercal~A\lr{~m+~S~y}\lr{~m+~S~y}^\intercal~B\lr{~m+~S~y} \\
        &= \lr{~y^\intercal~S~A~S~y}\lr{~y^\intercal~S~B~S~y} \\
        &\quad+ \lr{~m^\intercal~A~S~y}\lr{~y^\intercal~S~B~S~y} + \lr{~y^\intercal~S~A~m}\lr{~y^\intercal~S~B~S~y} + \lr{~y^\intercal~S~A~S~y}\lr{~m^\intercal~B~S~y} + \lr{~y^\intercal~S~A~S~y}\lr{~y^\intercal~S~B~m} \\
        &\quad + \lr{~m^\intercal~A~m}\lr{~y^\intercal~S~B~S~y} + \lr{~m^\intercal~A~S~y}\lr{~m^\intercal~B~S~y} + \lr{~m^\intercal~A~S~y}\lr{~y^\intercal~S~B~m} + \lr{~y^\intercal~S~A~m}\lr{~m^\intercal~B~S~y} \\
        &\quad\quad+ \lr{~y^\intercal~S~A~m}\lr{~y^\intercal~S~B~m} + \lr{~y^\intercal~S~A~S~y}\lr{~m^\intercal~B~m} \\
        &\quad+ \lr{~m^\intercal~A~m}\lr{~m^\intercal~B~S~y} + \lr{~m^\intercal~A~m}\lr{~y^\intercal~S~B~m} + \lr{~m^\intercal~A~S~y}\lr{~m^\intercal~B~m} + \lr{~y^\intercal~S~A~m}\lr{~m^\intercal~B~m} \\
        &\quad+ \lr{~m^\intercal~A~m}\lr{~m^\intercal~B~m}.
    \end{align*}
    Taking the expectation, it follows as in Lemma~\ref{lem:mean0quadprod} that odd-degree terms in $~y$ will drop out, therefore
    \begin{align*}
     \expec{\lr{~x^\intercal~A~x}\lr{~x^\intercal~B~x}} &= \lr{~m^\intercal~A~m}\lr{~m^\intercal~B~m} + \expec{\lr{~y^\intercal~S~A~S~y}\lr{~y^\intercal~S~B~S~y}} \\
     &\quad+ \lr{~m^\intercal~A~m}\expec{~y^\intercal~S~B~S~y} + \lr{~m^\intercal~B~m}\expec{~y^\intercal~S~A~S~y} + 4\,\expec{\lr{~m^\intercal~A~S~y}\lr{~m^\intercal~B~S~y}} \\
     &= \lr{~m^\intercal~A~m}\lr{~m^\intercal~B~m} + 2\IP{~S~A~S}{~S~B~S} + \tr\lr{~S~A~S}\tr\lr{~S~B~S} \\
     &\quad+ \lr{~m^\intercal~A~m}\tr\lr{~S~B~S} + \lr{~m^\intercal~B~m}\tr\lr{~S~A~S} + 4\,~m^\intercal~A~S^2~B~m \\
     &= \lr{~m^\intercal~A~m}\lr{~m^\intercal~B~m} + 2\IP{~A~C}{~C~B} + \tr\lr{~A~C}\tr\lr{~B~C} \\
     &\quad+ \lr{~m^\intercal~A~m}\tr\lr{~B~C} + \lr{~m^\intercal~B~m}\tr\lr{~A~C} + 4~m^\intercal~A~C~B~m,
    \end{align*}
where we have applied $~S^2=~C$ along with both Lemma~\ref{lem:quadform} and Lemma~\ref{lem:mean0quadprod}.  The result now follows from immediately from grouping terms and symmetrizing $~A,~B$.
\end{proof}

\begin{lemma}\label{lem:bigquadprod}
    Let $~x\sim\mathcal{N}\lr{~m, ~C}$, then for generic square $~A,~B$ and vectors $~a,~b,~c,~d$ it follows that 
    \begin{align*}
        \mathbb{E}\big[\lr{~x-~c}^\intercal~A&\lr{~x-~a}\lr{~x-~d}^\intercal~B\lr{~x-~b}\big] \\
        &= \IP{\lr{~A+~A^\intercal}~C}{~C~B} + \IP{~A}{~C}\IP{~B}{~C} + \lr{~m-~c}^\intercal~A\lr{~m-~a}\lr{~m-~d}^\intercal~B\lr{~m-~b} \\
        &\quad+ \lr{~m-~c}^\intercal~A\lr{~m-~a}\IP{~B}{~C} + \lr{~m-~d}^\intercal~B\lr{~m-~b}\IP{~A}{~C} \\
        &\quad+ \lr{~A\lr{~m-~a}+~A^\intercal\lr{~m-~c}}^\intercal~C\lr{~B\lr{~m-~b}+~B^\intercal\lr{~m-~d}}.
    \end{align*}
\end{lemma}
\begin{proof}
    First, notice that
    \begin{align*}
        \lr{~x-~c}^\intercal&~A\lr{~x-~a}\lr{~x-~d}^\intercal~B\lr{~x-~b} \\
        &= ~c^\intercal~A~a~d^\intercal~B~b - \lr{~x^\intercal~A~a~d^\intercal~B~b + ~c^\intercal~A~x~d^\intercal~B~b + ~c^\intercal~A~a~x^\intercal~B~b + ~c^\intercal~A~a~d^\intercal~B~x} \\
        &\quad+ ~x^\intercal~A~x~d^\intercal~B~b + ~x^\intercal~A~a~x^\intercal~B~b + ~x^\intercal~A~a~d^\intercal~B~x + ~c^\intercal~A~x~x^\intercal~B~b + ~c^\intercal~A~x~d^\intercal~B~x + ~c^\intercal~A~a~x^\intercal~B~x \\
        &\quad- \lr{~x^\intercal~A~x~x^\intercal~B~b + ~x^\intercal~A~x~d^\intercal~B~x + ~x^\intercal~A~a~x^\intercal~B~x + ~c^\intercal~A~x~x^\intercal~B~x} + ~x^\intercal~A~x~x^\intercal~B~x.
    \end{align*}
    We will calculate the expectation in stages, freely transposing scalar quantities and factoring whenever possible.  First, isolating the zero-order and first-order terms, we have
    \begin{align*}
        \mathbb{E}\big[ ~c^\intercal~A&~a~d^\intercal~B~b - \lr{~x^\intercal~A~a~d^\intercal~B~b + ~c^\intercal~A~x~d^\intercal~B~b + ~c^\intercal~A~a~x^\intercal~B~b + ~c^\intercal~A~a~d^\intercal~B~x} \big] \\
        &= ~c^\intercal~A~a~d^\intercal~B~b - \lr{\expec{~x}^\intercal~A~a~d^\intercal~B~b + ~c^\intercal~A\,\expec{~x}~d^\intercal~B~b + ~c^\intercal~A~a\,\expec{~x}^\intercal~B~b + ~c^\intercal~A~a~d^\intercal~B\,\expec{~x}} \\
        &= ~c^\intercal~A~a~d^\intercal~B~b - ~m^\intercal~A~a~d^\intercal~B~b - ~c^\intercal~A~m~d^\intercal~B~b - ~c^\intercal~A~a~m^\intercal~B~b - ~c^\intercal~A~a~d^\intercal~B~m.
    \end{align*}
    Now, isolating the second-order terms, it follows that
    \begin{align*}
        \mathbb{E}\big[~x^\intercal~A&~x~d^\intercal~B~b + ~x^\intercal~A~a~x^\intercal~B~b + ~x^\intercal~A~a~d^\intercal~B~x + ~c^\intercal~A~x~x^\intercal~B~b + ~c^\intercal~A~x~d^\intercal~B~x + ~c^\intercal~A~a~x^\intercal~B~x\big] \\
        &= \expec{~x^\intercal~A~x}~d^\intercal~B~b + ~a^\intercal~A^\intercal\,\expec{~x~x^\intercal}~B~b + ~d^\intercal~B\,\expec{~x~x^\intercal}~A~a + ~c^\intercal~A\,\expec{~x~x^\intercal}~B~b \\
        &\quad+ ~d^\intercal~B\,\expec{~x~x^\intercal}~A^\intercal~c + ~c^\intercal~A~a\,\expec{~x^\intercal~B~x} \\
        &= ~m^\intercal~A~m~d^\intercal~B~b + ~m^\intercal~A~a~m^\intercal~B~b+~m^\intercal~A~a~d^\intercal~B~m+~c^\intercal~A~m~m^\intercal~B~b+~c^\intercal~A~m~d^\intercal~B~m + ~c^\intercal~A~a~m^\intercal~B~m \\
        &\quad+ \IP{~A}{~C}~d^\intercal~B~b + \IP{~B}{~C}~c^\intercal~A~a + \lr{~a^\intercal~A^\intercal + ~c^\intercal~A}~C\lr{~B~b + ~B^\intercal~d}.
    \end{align*}
    Similarly, isolating the third-order and fourth-order terms gives
    \begin{align*}
        \mathbb{E}\big[~x^\intercal~A&~x~x^\intercal~B~x - ~x^\intercal~A~x~x^\intercal~B~b - ~x^\intercal~A~x~d^\intercal~B~x - ~x^\intercal~A~a~x^\intercal~B~x - ~c^\intercal~A~x~x^\intercal~B~x. \big] \\
        &= \expec{~x^\intercal~A~x~x^\intercal~B~x} - \expec{~x^\intercal~A~x~x^\intercal}~B~b - ~d^\intercal~B\,\expec{~x~x^\intercal~A~x} - \expec{~x^\intercal~B~x~x^\intercal}~A~a - ~c^\intercal~A\,\expec{~x~x^\intercal~B~x} \\
        &= \IP{\lr{~A+~A^\intercal}~C}{~C~B} + ~m^\intercal\lr{~A+~A^\intercal}~C\lr{~B+~B^\intercal}~m + \lr{~m^\intercal~A~m + \IP{~A}{~C}}\lr{~m^\intercal~B~m + \IP{~B}{~C}} \\
        &\quad- \lr{\IP{~A}{~C}~m^\intercal~B~b + ~m^\intercal\lr{~A+~A^\intercal}~C~B~b + ~m^\intercal~A~m~m^\intercal~B~b} \\
        &\quad- \lr{~d^\intercal~B~m\IP{~A}{~C} + ~d^\intercal~B~C\lr{~A+~A^\intercal}~m + ~d^\intercal~B~m~m^\intercal~A~m} \\
        &\quad- \lr{\IP{~B}{~C}~m^\intercal~A~a + ~m^\intercal\lr{~B+~B^\intercal}~C~A~a + ~m^\intercal~B~m~m^\intercal~A~a} \\
        &\quad- \lr{~c^\intercal~A~m\IP{~B}{~C} + ~c^\intercal~A~C\lr{~B+~B^\intercal}~m + ~c^\intercal~A~m~m^\intercal~B~m} \\
        &= \IP{\lr{~A+~A^\intercal}~C}{~C~B} + ~m^\intercal\lr{~A+~A^\intercal}~C\lr{~B+~B^\intercal}~m + \IP{~A}{~C}\IP{~B}{~C} \\
        &\quad+ \IP{~A}{~C}\lr{~m^\intercal~B~m-~m^\intercal~B~b-~d^\intercal~B~m} + \IP{~B}{~C}\lr{~m^\intercal~A~m-~m^\intercal~A~a-~c^\intercal~A~m} \\
        &\quad- ~m^\intercal\lr{~A+~A^\intercal}~C\lr{~B~b+~B^\intercal~d} - \lr{~a^\intercal~A^\intercal+~c^\intercal~A}~C\lr{~B+~B^\intercal}~m \\
        &\quad+ ~m^\intercal~A~m~m^\intercal~B~m - ~m^\intercal~A~m~m^\intercal~B~b - ~m^\intercal~A~m~d^\intercal~B~m - ~m^\intercal~A~a~m^\intercal~B~m - ~c^\intercal~A~m~m^\intercal~B~m
    \end{align*}
    Putting this all together, it follows that the expectation is computable as 
    \begin{align*}
        \mathbb{E}\big[ (~x&-~c)^\intercal~A\lr{~x-~a}\lr{~x-~d}^\intercal~B\lr{~x-~b}\big] \\
        &= \IP{\lr{~A+~A^\intercal}~C}{~C~B} + \IP{~A}{~C}\IP{~B}{~C} \\
        &\quad+ \big( ~m^\intercal~A~m~m^\intercal~B~m - ~m^\intercal~A~m~m^\intercal~B~b - ~m^\intercal~A~m~d^\intercal~B~m + ~m^\intercal~A~m~d^\intercal~B~b \\
        &\quad\quad- ~m^\intercal~A~a~m^\intercal~B~m + ~m^\intercal~A~a~m^\intercal~B~b + ~m^\intercal~A~a~d^\intercal~B~m -~m^\intercal~A~a~d^\intercal~B~b \\
        &\quad\quad- ~c^\intercal~A~m~m^\intercal~B~m + ~c^\intercal~A~m~m^\intercal~B~b + ~c^\intercal~A~m~d^\intercal~B~m - ~c^\intercal~A~m~d^\intercal~B~b \\
        &\quad\quad+ ~c^\intercal~A~a~m^\intercal~B~m - ~c^\intercal~A~a~m^\intercal~B~b - ~c^\intercal~A~a~d^\intercal~B~m + ~c^\intercal~A~a~d^\intercal~B~b \big) \\
        &\quad+ \IP{~B}{~C}\lr{~c^\intercal~A~a - ~m^\intercal~A~a - ~c^\intercal~A~m + ~m^\intercal~A~m} + \IP{~A}{~C}\lr{~d^\intercal~B~b -~m^\intercal~B~b - ~d^\intercal~B~m + ~m^\intercal~B~m} \\
        &\quad+ \big( \lr{~a^\intercal~A^\intercal + ~c^\intercal~A}~C\lr{~B~b+~B^\intercal~d} - ~m^\intercal\lr{~A+~A^\intercal}~C\lr{~B~b+~B^\intercal~d} \\
        &\quad\quad- \lr{~a^\intercal~A^\intercal + ~c^\intercal~A}~C\lr{~B+~B^\intercal}~m + ~m^\intercal\lr{~A+~A^\intercal}~C\lr{~B+~B^\intercal}~m \big).
    \end{align*}
    Examining each group of terms finally yields the conclusion. 
\end{proof}

\subsubsection{Two-form mass matrix}
\label{subsec:two-forms}
Recall that if Gaussian densities $\phi_i(~x)$ are chosen as zero-forms, then we have $\nabla\phi_i(~x) = -\phi_i(~x)~C_i^{-1}\lr{~x-~m_i}$.  It follows that (a proxy for) the basis Whitney 2-form is given by
\begin{align*}
    \psi_{ijk} &= \phi_i\nabla\phi_j\times\nabla\phi_k + \phi_j\nabla\phi_k\times\nabla\phi_i + \phi_k\nabla\phi_i\times\nabla\phi_j \\
    &= \phi_i\phi_j\phi_k \sum_{(ijk)} ~C_i^{-1}\lr{~x-~m_i}\times~C_j^{-1}\lr{~x-~m_j},
\end{align*}
where the sum over $(ijk)$ denotes a sum over cyclic permutations of these indices.  The goal is to compute the mass matrix
\begin{align*}
    \mathcal{M}_{ijk,abc} &= \int_M \psi_{ijk}\cdot\psi_{abc} \\
    &= \int_M \phi_i\phi_j\phi_k\phi_a\phi_b\phi_c\sum_{(ijk)}~C_i^{-1}\lr{~x-~m_i}\times~C_j^{-1}\lr{~x-~m_j}\cdot\sum_{(abc)}~C_a^{-1}\lr{~x-~m_a}\times~C_b^{-1}\lr{~x-~m_b} \\
    &= z_{ijkabc}\int_M \phi_{ijkabc}(~x)\sum_{(ijk)}\sum_{(abc)}\lr{~x-~m_i}^\intercal~C_i^{-1}~C_a^{-1}\lr{~x-~m_a}\lr{~x-~m_j}^\intercal~C_j^{-1}~C_b^{-1}\lr{~x-~m_b} \\
    &\quad- z_{ijkabc}\int_M \phi_{ijkabc}(~x)\sum_{(ijk)}\sum_{(abc)}\lr{~x-~m_j}^\intercal~C_j^{-1}~C_a^{-1}\lr{~x-~m_a}\lr{~x-~m_i}^\intercal~C_i^{-1}~C_b^{-1}\lr{~x-~m_b} \\
    &= z_{ijkabc}\sum_{(ijk)}\sum_{(abc)}\expec{\lr{~x-~m_i}^\intercal~C_i^{-1}~C_a^{-1}\lr{~x-~m_a}\lr{~x-~m_j}^\intercal~C_j^{-1}~C_b^{-1}\lr{~x-~m_b}} \\
    &\quad- z_{ijkabc}\sum_{(ijk)}\sum_{(abc)}\expec{\lr{~x-~m_j}^\intercal~C_j^{-1}~C_a^{-1}\lr{~x-~m_a}\lr{~x-~m_i}^\intercal~C_i^{-1}~C_b^{-1}\lr{~x-~m_b}},
\end{align*}
where, again, the expectation is taken with respect to $\phi_{ijkabc}(~x)$.  The normalization term $z_{ijkabc}$, mean $~m$, and covariance $~C$ associated to the Gaussian density $\phi_{ijkabc}$ have already been computed in Lemma~\ref{lem:gaussprod}.  Applying Lemma~\ref{lem:bigquadprod}, it now follows that
\begin{align*}
    \mathbb{E}\big[ (~x&-~m_i)~C_i^{-1}~C_a^{-1}\lr{~x-~m_a}\lr{~x-~m_j}^\intercal~C_j^{-1}~C_b^{-1}\lr{~x-~m_b} \big] \\
    &= \IP{\lr{~C_i^{-1}~C_a^{-1}+~C_a^{-1}~C_i^{-1}}~C}{~C~C_j^{-1}~C_b^{-1}} + \IP{~C_i^{-1}~C_a^{-1}}{~C}\IP{~C_j^{-1}~C_b^{-1}}{~C} \\
    &\quad+ \lr{~m-~m_i}~C_i^{-1}~C_a^{-1}\lr{~m-~m_a}\lr{~m-~m_j}^\intercal~C_j^{-1}~C_b^{-1}\lr{~m-~m_b} \\
    &\quad+ \lr{~m-~m_i}^\intercal~C_i^{-1}~C_a^{-1}\lr{~m-~m_a}\IP{~C_j^{-1}~C_b^{-1}}{~C} \\
    &\quad+ \lr{~m-~m_j}^\intercal~C_j^{-1}~C_b^{-1}\lr{~m-~m_b}\IP{~C_i^{-1}~C_a^{-1}}{~C} \\
    &\quad+ \lr{~C_i^{-1}~C_a^{-1}\lr{~m-~m_a}+~C_a^{-1}~C_i^{-1}\lr{~m-~m_i}}^\intercal~C\lr{~C_j^{-1}~C_b^{-1}\lr{~m-~m_b}+~C_b^{-1}~C_j^{-1}\lr{~m-~m_j}}.
\end{align*}
Similarly, the second term is 
\begin{align*}
    \mathbb{E}\big[ (~x&-~m_j)~C_j^{-1}~C_a^{-1}\lr{~x-~m_a}\lr{~x-~m_i}^\intercal~C_i^{-1}~C_b^{-1}\lr{~x-~m_b} \big] \\
    &= \IP{\lr{~C_j^{-1}~C_a^{-1}+~C_a^{-1}~C_j^{-1}}~C}{~C~C_i^{-1}~C_b^{-1}} + \IP{~C_j^{-1}~C_a^{-1}}{~C}\IP{~C_i^{-1}~C_b^{-1}}{~C} \\
    &\quad+ \lr{~m-~m_j}~C_j^{-1}~C_a^{-1}\lr{~m-~m_a}\lr{~m-~m_i}^\intercal~C_i^{-1}~C_b^{-1}\lr{~m-~m_b} \\
    &\quad+ \lr{~m-~m_j}^\intercal~C_j^{-1}~C_a^{-1}\lr{~m-~m_a}\IP{~C_i^{-1}~C_b^{-1}}{~C} \\
    &\quad+ \lr{~m-~m_i}^\intercal~C_i^{-1}~C_b^{-1}\lr{~m-~m_b}\IP{~C_j^{-1}~C_a^{-1}}{~C} \\
    &\quad+ \lr{~C_j^{-1}~C_a^{-1}\lr{~m-~m_a}+~C_a^{-1}~C_j^{-1}\lr{~m-~m_j}}^\intercal~C\lr{~C_i^{-1}~C_b^{-1}\lr{~m-~m_b}+~C_b^{-1}~C_i^{-1}\lr{~m-~m_i}}.
\end{align*}
Denoting these terms as $T_{iajb}$ and $T_{jaib}$, respectively, it follows that the mass matrix is computable through
\[~\mathcal{M}_{ijk,abc} = \int_M \psi_{ijk}\cdot\psi_{abc} = z_{ijkabc}\sum_{(ijk)}\sum_{(abc)} \lr{T_{iajb} - T_{jaib}},\]
where $z_{ijkabc}, ~C = ~C_{ijkabc},$ and $~m=~m_{ijkabc}$ are computable via Lemma~\ref{lem:gaussprod}.

\end{document}